\documentclass[a4paper,fleqn]{cas-sc}

\usepackage[authoryear]{natbib}

\def\tsc#1{\csdef{#1}{\textsc{\lowercase{#1}}\xspace}}
\tsc{WGM}
\tsc{QE}

\usepackage{graphicx}
\usepackage{times}
\usepackage{soul}
\usepackage{url}
\usepackage{booktabs}
\usepackage{placeins}
\usepackage{algorithm}
\usepackage{algorithmic}
\usepackage{amssymb}
\usepackage{subcaption}
\usepackage{lineno}
\usepackage{enumitem}
\usepackage{caption}
\usepackage{adjustbox}

\newtheorem{proposition}{Proposition}
\newtheorem{lemma}{Lemma}

\DeclareMathOperator*{\argmin}{arg\,min}
\newcommand{\E}{\mathbb{E}}

\begin{document}
\let\WriteBookmarks\relax
\def\floatpagepagefraction{1}
\def\textpagefraction{.001}

\shorttitle{CPDs for Contextual Order Fulfillment}

\shortauthors{}  

\renewcommand{\printorcid}{}
\ExplSyntaxOn
\cs_set:Npn \__first_footerline: { \group_begin: \small \sffamily \__short_authors: \group_end: }
\ExplSyntaxOff
\shorttitle{Conformal Predictive Distributions for Order Fulfillment}
\shortauthors{T. Ye et al.}
\title [mode = title]{Conformal Predictive Distributions for Contextual Order Fulfillment Optimization}
\tnotemark[1]
\tnotetext[1]{A conference version of this paper appeared as T. Ye, A. Hijazi, and P. Van Hentenryck, ``Conformal Predictive Distributions for Order Fulfillment Time Forecasting,'' Joint International Conference on Computational Logistics and EURO Mini Conference on Maritime Optimization and Logistics (ICCL--EUROMar 2025), Springer Lecture Notes in Computer Science. The forecasting study of Section \ref{sec:experiments} first appeared in that version.}

\author[1]{Tinghan Ye}
\cormark[1]
\ead{joe.ye@gatech.edu}

\author[1]{Amira Hijazi}
\ead{ahijazi6@gatech.edu}

\author[1]{Pascal {Van Hentenryck}}
\ead{pvh@gatech.edu}

\affiliation[1]{organization={H. Milton Stewart School of Industrial \& Systems Engineering, Georgia Institute of Technology},
                city={Atlanta},
                state={Georgia},
                country={USA}}

\cortext[1]{Corresponding author}

\begin{abstract}
Ship-to-home fulfillment requires a retailer to select a fulfillment location and carrier before the delivery time is known. These decisions often rely on static transit-time estimates, which cannot represent how delivery-time uncertainty varies across orders and fulfillment options. This paper constructs conformal predictive distributions (CPDs) using Conformal Predictive Systems (CPS) for regression and the cross Venn--Abers predictor (CVAP) for multiclass classification, and uses them in a contextual optimization model for fulfillment decisions. The model is solved by contextual sample-average approximation (C-SAA). The regret of each fulfillment decision, its expected cost above that of the best decision, is bounded by forecast error, scenario-evaluation error, and the candidate-generation gap. The forecast term depends on how well calibration residuals or class probabilities represent the current order. On more than five million proprietary order records, CVAP raises exact-day accuracy from 49.5\% for the rule-based estimate to 64.1\%. Paired with a decision rule that selects a quantile weighting late errors more heavily, the split CPS (SCPS) gives a point estimate that raises late-delivery detection from 3.8\% to 79.4\%. Residual-based CPS distributions also meet or slightly exceed nominal interval-coverage targets. In decision experiments with a JD.com simulator over six test dates, C-SAA with CVAP lowers mean total cost by 7.1\% and C-SAA with the Mondrian CPS (MCPS) lowers it by 6.7\%, relative to C-SAA with the original nonconformal forecast. Both cut the late-delivery rate from 5.7\% to below 4.0\%.
\end{abstract}


\maketitle

\section{Introduction}
\begin{figure*}[pos=!ht]
    \centering \includegraphics[width=\textwidth]{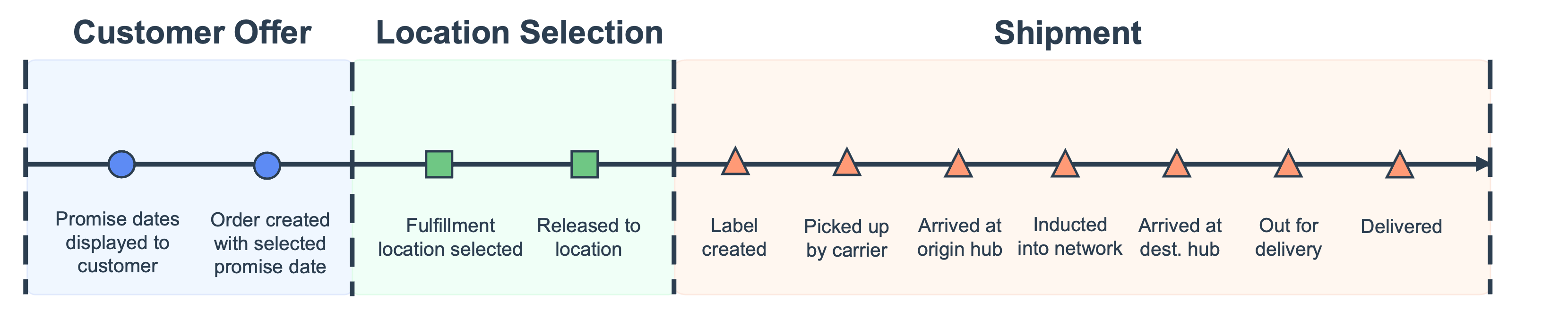}
    \caption{Example of Order Fulfillment Life Cycle.}
    \label{fig:life-cycle}
\end{figure*}

Ship-to-home fulfillment requires a retailer to choose a fulfillment option, a pairing of a fulfillment location with a carrier service, before the delivery outcome is known. The fulfillment decision determines shipping expense, opportunities to consolidate items, and the likelihood that the order arrives early or late. Many systems make this choice using static transit tables that specify service times, pickup cutoffs, and valid operating days. These tables are easy to maintain, but they do not adapt to the order, current fulfillment operations, or conditions in the carrier network. Machine-learning methods are increasingly used to predict delivery times, but accurate prediction remains difficult because the outcome depends on every stage the order passes through: processing and shipping at the fulfillment location, entry into the carrier network, line-haul transport, and the last-mile leg to the customer \citep{salari2022real,merchan20242021,ruan2022service}. Figure \ref{fig:life-cycle} traces these stages from order placement through delivery.

For the fulfillment decision, the form of the forecast matters. A point forecast supplies one estimated delivery date but not the likelihood or extent of an early or late arrival. The difference is important when the two types of deviation carry asymmetric costs and when fulfillment options with similar expected delivery dates have different levels of uncertainty. A predictive distribution supplies the scenarios needed by contextual stochastic optimization and can also be summarized as a point forecast, quantile, or interval \citep{vovk2020computationally}. The operational question is how the representation of delivery-time uncertainty affects fulfillment decisions.

Two related studies motivate this question. \citet{ye2025msom} formulate contextual stochastic optimization for omnichannel multicourier fulfillment and use estimated delivery-time distributions to account for item consolidation and delivery uncertainty. They estimate these distributions with a multiclass classifier for each carrier, calibrated by isotonic regression, and with tree-based quantile regression, neither of which carries a finite-sample validity guarantee. \citet{ye2025iccl} develop conformal predictive distributions (CPDs) for fulfillment-time forecasting and evaluate their distributional and point accuracy on proprietary order records. The present study uses those forecasting methods to construct inputs for the stochastic fulfillment model, bounds the cost of the resulting fulfillment decisions, and evaluates the policies they produce in the simulator of \citet{ye2026proxy}.

The forecasting methods of \citet{ye2025iccl} construct CPDs in two ways. Conformal Predictive Systems (CPS) convert regression residuals, the differences between observed outcomes and fitted predictions, into a predictive distribution. Split CPS (SCPS) uses all calibration residuals, whereas Mondrian CPS (MCPS) uses residuals from groups with similar predicted values. The cross Venn--Abers predictor (CVAP) constructs a CPD by calibrating probabilities from a multiclass delivery-deviation model. The forecasting comparison also includes two-stage classify-then-regress forecasts and delivery-date estimates chosen to reflect the different consequences of early and late errors. In the decision study, the CPDs supply delivery scenarios to a contextual fulfillment model solved by contextual sample-average approximation (C-SAA). C-SAA evaluates feasible fulfillment decisions by averaging their costs over sampled delivery and future-demand outcomes. Using the same solution procedure for every forecast focuses the downstream comparison on the distribution supplied to the optimizer.

The methodological contribution is to represent delivery-time uncertainty with CPDs for each order and fulfillment option, the level at which the optimization model compares fulfillment options. A general result bounds the regret of the implemented decision, the amount by which its expected cost exceeds that of the best decision for the current order. The bound combines the error in the expected delivery cost that the forecast assigns to each fulfillment option, the sampling error in evaluating candidate decisions, and the loss from choosing among a finite set of candidates. Two further results bound the error in expected delivery cost. For CPS, the bound depends on calibration-group size, regression error, and how well calibration residuals match the current order. For class-probability forecasts such as CVAP, it depends on their accuracy at the current order. An example shows that exact marginal calibration does not bound this error, so these results identify the order-level accuracy that the decision depends on. The empirical contribution evaluates these forecasts on more than five million order records and compares the fulfillment policies they produce in a simulator.

\section{Related Work}
\label{sec:related-work}

This paper connects three literatures. Delivery-time forecasting supplies the inputs, probabilistic and conformal forecasting determine their form, and contextual optimization turns them into fulfillment decisions.

\subsection{Delivery-Time Forecasting}
Delivery-time forecasting has been studied at the customer-offer, fulfillment, line-haul, and last-mile stages. \citet{wen2024survey} provide a recent review. Methods include deep travel-time models \citep{wu2019deepeta,de2021end,li2021unsupervised,ruan2022service}, tree-based models \citep{salari2022real,zhang2023delivery}, and models tailored to delivery exceptions and service operations \citep{cao2023unraveling,yi2024learning,zhang2024estimating}. Most of this literature predicts the customer promise date. This paper instead predicts delivery-time deviation when a fulfillment option is considered, before the fulfillment decision is made. The timing of the forecast follows \citet{ye2025msom}, but the present study replaces its incumbent distributional forecasts with CPDs and adds carrier-network features.

\subsection{Probabilistic and Conformal Forecasting}
Probabilistic forecasting represents delivery-time uncertainty through a predictive distribution. Parametric service-time models specify a distributional family for the outcome \citep{shang2017exploiting,liu2023uncertainty,raj2024stochastic}. Other approaches, including quantile forests \citep{meinshausen2006quantile,guo2022forecasting,salari2022real,arora2023probabilistic} and neural distributional models \citep{liu2023uncertainty,faulkner2025framework}, estimate the distribution more flexibly. The resulting distributions take two forms. Regression-based methods often estimate a set of delivery-time quantiles, such as the median and upper-tail percentiles, whereas multiclass models assign a probability to each possible integer-day deviation.

Conformal prediction adds a finite-sample guarantee to such forecasts. Under exchangeability, it provides marginal calibration, a guarantee averaged over calibration and test examples, without assuming a parametric outcome distribution \citep{vovk2005algorithmic,manokhin_valery_2022_6467205}. For regression, CPS ranks possible outcomes against held-out calibration residuals to construct a predictive cumulative distribution function (CDF) with this finite-sample guarantee \citep{vovk2020computationally}. Venn--Abers was developed for binary classification and returns a pair of probabilities with a finite-sample validity guarantee \citep{vovk2015large,manokhin2017multi}. Related logistics applications use conformalized quantile regression in a classify-then-regress model \citep{faulkner2025framework} or construct prediction sets for inbound load planning \citep{bruys2024confidence}.

\subsection{Contextual Optimization}
Contextual optimization uses information available before a decision to model uncertainty and choose an action. Feature-based newsvendor models and general predictive-to-prescriptive methods helped establish this literature \citep{ban2019big,bertsimas2020predictive}. \citet{sadana2025survey} review its main approaches and applications. The approach closest to the CPS construction is residual-based contextual sample-average approximation (SAA): a regression model supplies a prediction for the current context, and shifted residuals provide scenarios for the stochastic program. \citet{kannan2022data} analyze three such data-driven SAA constructions and establish consistency, convergence rates, and finite-sample guarantees under conditions on the data-generating process, prediction model, and stochastic program. \citet{sun2023predict} calibrate held-out residuals of an arbitrary prediction model for robust contextual linear programs. Their coverage guarantees are marginal, and their expected-cost rate assumes that the conditional mean residual varies smoothly with the covariates.

CPS resembles residual-based contextual SAA because it also adds calibration residuals to a regression prediction, but the two methods attach different guarantees to this construction. Residual-based SAA uses residuals to approximate the conditional distribution needed by the optimization model, whereas CPS ranks calibration residuals to obtain marginal predictive validity. This validity is marginal because distribution-free guarantees conditional on the covariates require additional structure \citep{barber2021limits}.

Other studies use conformal methods to guide decisions directly. Under finite decisions and bounded utility, \citet{vovk2018decision} provide an asymptotic decision result for universally consistent predictive systems. Another line of work links conformal prediction to robust optimization by using a calibrated prediction region as the uncertainty set \citep{johnstone2021conformal,patel2024conformal,yeh2025endtoend}, or chooses actions for risk-averse agents from conformal sets \citep{kiyani2025decision,zhu2026action}. Those approaches protect a decision against outcomes in a region, with guarantees that hold marginally or conditionally on the chosen action. By contrast, the model here samples full predictive distributions within a stochastic program to compare expected fulfillment costs.

A related literature asks what calibration implies for decisions. A forecast that ignores the covariates can be perfectly calibrated \citep{gneiting2007probabilistic}, and average notions of calibration need not give accurate estimates of decision loss \citep{zhao2021calibrating,sahoo2021reliable}. These studies propose calibration notions tied to classes of decisions, but their guarantees remain averages over the population. By contrast, the results in Section \ref{sec:theory-guarantees} bound the error in expected cost at the fulfillment options of the current order.

\section{Problem Setting}
\label{sec:problem-setting}

When an order arrives, the retailer chooses among fulfillment options before the delivery time is known. Each fulfillment option pairs a fulfillment location with a carrier service and differs in shipping cost, delivery-time uncertainty, and use of inventory. The problem therefore has two linked components: a delivery-time forecast for each fulfillment option and a sequential fulfillment decision that accounts for current operating conditions and later orders.

\subsection{Fulfillment Problem}

Orders arrive sequentially over a planning horizon of $T$ decision epochs, indexed by $t=1,\ldots,T$. Let $\mathcal S$ be the state space and $\boldsymbol s_t\in\mathcal S$ the system state when order $t$ arrives. The state describes the current order by its requested items and quantities, destination, arrival time, and promised delivery date. It also records the available inventory, the eligible carrier services, and the workload of each fulfillment location, measured by its recently shipped and pending units. Finally, the state contains the order, location, and carrier predictors described in Section \ref{sec:forecasting-problem}. The forecasting models receive these predictors, whereas the fulfillment decision uses the entire state. The set $\mathcal A(\boldsymbol s_t)$ contains the eligible fulfillment options, and $\mathcal Z(\boldsymbol s_t)$ contains the feasible fulfillment decisions. Before delivery outcomes are known, the retailer chooses $\boldsymbol z_t\in\mathcal Z(\boldsymbol s_t)$. This decision assigns each requested unit to an eligible fulfillment option or leaves it unfulfilled, and the units of an item assigned to a fulfillment location cannot exceed its available inventory there. Because the fulfillment decision consumes inventory, it affects the resources available for later orders.

For each fulfillment option $a\in\mathcal A(\boldsymbol s_t)$, let $Y_{a,t}$ denote the actual delivery date minus the promised delivery date, measured in days, if that fulfillment option is used. A negative value denotes an early delivery, zero an on-time delivery, and a positive value a late delivery. The vector $\boldsymbol Y_t=(Y_{a,t})_{a\in\mathcal A(\boldsymbol s_t)}$ collects these delivery-time deviations. The current-order cost $c(\boldsymbol s_t,\boldsymbol z_t,\boldsymbol Y_t)$ has two parts. An order-level term charges unfulfilled units and credits consolidation discounts. Each fulfillment option used adds a delivery cost, covering shipping and early- and late-delivery penalties, that depends on its deviation $Y_{a,t}$.

\subsection{Delivery-Time Forecasting Problem}
\label{sec:forecasting-problem}

Given the state $\boldsymbol s_t$, the delivery-time deviations $\boldsymbol Y_t$ of the eligible fulfillment options have a joint conditional distribution. Following \citet{ye2025msom}, this joint distribution is decomposed into one distribution per fulfillment option. Let the feature vector $\boldsymbol x_{a,t}$ be the part of $\boldsymbol s_t$ used to forecast the delivery time of fulfillment option $a$, consisting of predictors for the current order, the option's fulfillment location, and its carrier service. Given the state, the deviations are independent across fulfillment options, and each depends on the state only through its own feature vector:
\[
\mathbb P(\boldsymbol Y_t\leq\boldsymbol y\mid\boldsymbol s_t)=\prod_{a\in\mathcal A(\boldsymbol s_t)}F(y_a\mid\boldsymbol x_{a,t}),\qquad \boldsymbol y=(y_a)_{a\in\mathcal A(\boldsymbol s_t)},
\]
where $F(y\mid\boldsymbol x_{a,t})=\mathbb P(Y_{a,t}\leq y\mid\boldsymbol x_{a,t})$ is the conditional CDF of the deviation at fulfillment option $a$. Because the current-order cost adds up the delivery costs of the fulfillment options used, these per-option distributions determine the expected current-order cost.

The forecasting target is therefore $F(\cdot\mid\boldsymbol x_{a,t})$, and each method returns an estimate $\widehat F(\cdot\mid\boldsymbol x_{a,t})$. Two kinds of models produce this estimate. Regression models predict either the deviation $Y_{a,t}$ or the delivery duration, and a duration forecast is converted to a deviation from the promised date before delivery penalties are applied. Let $K$ be the number of distinct integer-day outcomes observed in training. The multiclass classifier (MC-CLF) instead treats these $K$ outcomes as classes and assigns a probability to each, and these probabilities define a discrete CDF. The estimate $\widehat F$ can be sampled to generate optimization scenarios or summarized by its quantiles: one quantile gives a point estimate and two give a prediction interval.

\subsection{Sequential Contextual Stochastic Formulation}
\label{sec:decision-problem}

Let $O_{t+1}$ denote the next order, with $O_{T+1}$ empty because the horizon ends after order $T$. After decision $\boldsymbol z_t$, the state moves to $\boldsymbol s_{t+1}=f(\boldsymbol s_t,\boldsymbol z_t,O_{t+1})$, where the transition function $f$ removes the assigned units from inventory, updates the workload of the chosen fulfillment locations, and records the next order. The delivery-time deviations $\boldsymbol Y_t$ affect the current-order cost but not the next state.

A policy $\pi=(\pi_1,\ldots,\pi_T)$ specifies one decision rule for each epoch, where $\pi_t$ maps each state $\boldsymbol s_t\in\mathcal S$ to a feasible decision $\pi_t(\boldsymbol s_t)\in\mathcal Z(\boldsymbol s_t)$. The sequential contextual stochastic problem is
\begin{equation}
\min_{\pi}\ \E\!\left[\sum_{t=1}^{T}
c\bigl(\boldsymbol s_t,\pi_t(\boldsymbol s_t),\boldsymbol Y_t\bigr)\right],
\label{eq:sequential-problem}
\end{equation}
where the initial state $\boldsymbol s_1$ is given. The expectation is over future orders, which determine the states visited under $\pi$, and over the delivery-time deviations, with each $Y_{a,t}$ following $F(\cdot\mid\boldsymbol x_{a,t})$.

\FloatBarrier

\section{Delivery-Time Forecasting Methods}
\label{sec:forecast-methodology}
This section describes how the forecasts estimate the conditional distribution $F(\cdot\mid\boldsymbol x_{a,t})$ of Section \ref{sec:problem-setting}. The estimate assigns probabilities across early, on-time, and late outcomes, so it can be used to construct prediction intervals, sample delivery scenarios for contextual optimization, or provide a single estimated date. Regression and multiclass classification give two routes to such a distribution.

\subsection{Conformal Predictive Distributions for Regression and Classification} \label{sec:conformal}
Regression and classification models produce different outputs, so their calibration steps differ. CPS turns a regression forecast and its past errors into a predictive distribution. Venn--Abers calibration adjusts the probabilities assigned by MC-CLF to the possible delivery-day classes. Both approaches express the result as an estimated CDF over delivery-time deviations.

\subsubsection{Conformal Predictive System} \label{sec:cps}

The \textit{Conformal Predictive System for Regression} (CPS) starts with a regression estimate and describes uncertainty around it using errors on held-out calibration orders. Adding those errors to the estimate for a new order gives the delivery times that would follow if calibration-period errors recurred, and the resulting predictive distribution supplies quantiles for point forecasts and prediction intervals.

Let $h$ be the regression model trained on the proper training set, and let $N_c$ be the number of held-out calibration observations, indexed as $(\boldsymbol x_j,y_j)$ for $j=1,\ldots,N_c$. For a test input $\boldsymbol x$, CPS adds each calibration residual $y_j-h(\boldsymbol x_j)$ to the regression prediction $h(\boldsymbol x)$ to obtain a possible deviation $A_j(\boldsymbol x)$:
\(
A_j(\boldsymbol x)=h(\boldsymbol x)+y_j-h(\boldsymbol x_j).
\)
The deterministic CPS of \citet{vovk2020computationally} assigns probability $1/N_c$ to each shifted outcome, giving the predictive CDF
\begin{equation}
\widehat F(y\mid \boldsymbol x)
=\frac{1}{N_c}\sum_{j=1}^{N_c}\mathbf 1\{A_j(\boldsymbol x)\leq y\},
\label{eq:cps-cdf}
\end{equation}
where $\mathbf 1\{\cdot\}$ is the indicator function. C-SAA samples this distribution. \citet{vovk2020computationally} also define a randomized CPS, which differs from Equation \eqref{eq:cps-cdf} in how it treats ties and the candidate outcome itself. The randomized CPS is valid in finite samples. When the calibration and test examples are exchangeable, its value at the observed test outcome is uniformly distributed on $[0,1]$ (Appendix \ref{app:forecast-methods}).

SCPS uses all $N_c$ calibration errors in Equation \eqref{eq:cps-cdf} and therefore applies the same empirical error distribution to every new prediction. MCPS groups the calibration observations by ranges of $h(\boldsymbol x_j)$. For a test prediction in one range, MCPS evaluates Equation \eqref{eq:cps-cdf} with the calibration errors from that range \citep{bostrom2021mondrian}. Grouping allows MCPS to represent different error patterns for near-on-time predictions and predicted delays.

\subsubsection{Venn--Abers Calibration for Multiclass Classification (MC-CLF)} \label{sec:va}

MC-CLF assigns a probability to each of the $K$ deviation classes, and these probabilities determine the predictive CDF. Because the raw probabilities may be systematically too high or too low, calibration compares the model scores with observed class frequencies and adjusts the probabilities while preserving their order. Venn--Abers predictors, originally developed for binary classification, perform this adjustment by combining isotonic regression with conformal calibration \citep{vovk2015large,manokhin2017multi}. To calibrate the multiclass probabilities, this paper adopts a one-versus-rest strategy, applying a separate Venn--Abers predictor to each class against all remaining classes. The inductive Venn--Abers predictor (IVAP) considers both possible labels for a new instance. It appends the instance to the calibration data once with each label, fits a nondecreasing map from model scores to observed class frequencies each time, and thus obtains a pair of calibrated probabilities (Appendix \ref{app:forecast-methods}). The distance between these two values records how much the calibrated probability depends on the unknown label. Under exchangeability, the probability computed with the true label appended is perfectly calibrated \citep{vovk2015large}.

CVAP repeats IVAP over several training and calibration splits, called folds, to reduce dependence on a single split. Because the predictive CDF needs one probability, CVAP merges the fold-specific pairs into a single forecast with a geometric-mean rule. The merged forecast gives up the validity guarantee of IVAP, although merged Venn--Abers forecasts were more accurate than existing calibration methods in the experiments of \citet{vovk2015large}. The resulting one-versus-rest probabilities are then normalized to sum to one across the $K$ classes.

\subsection{2-Stage Classify-then-Regress Forecast}
\label{sec:classify-regress-forecast}

Delivery-time deviation has an ordered numerical structure: both its sign and its distance from zero are meaningful. The 2-stage forecast represents this structure by separating delivery status from deviation magnitude. A Venn--Abers-calibrated classifier first estimates the probabilities of early, on-time, and late delivery. The forecast then mixes three component distributions with these probabilities as weights: an early-delivery component retaining the negative outcomes of a CPS distribution, an on-time component concentrated at zero, and a late-delivery component retaining the positive outcomes. Using SCPS or MCPS for the early and late deviation distributions yields the 2-stage SCPS and 2-stage MCPS variants. Here \textit{2-stage} refers to the forecasting construction, not to two-stage stochastic optimization.

\subsection{Obtaining Point Predictions from Distributional Forecasts}
\label{sec:point-prediction}

A retailer can select a delivery-date estimate by choosing a quantile of the predictive distribution. The median targets its center, whereas a higher quantile shifts the estimated date later and places more emphasis on avoiding underestimation of delivery delays. The quantile can therefore be chosen to reflect the relative consequences of errors on early and late deliveries without refitting the forecasting model.

For CPS, the conventional point forecast is the median of the predicted CDF. For MC-CLF, it is the most probable deviation class. A weighted-error rule instead chooses a quantile level by balancing overall prediction error against errors on late and early deliveries. For a candidate quantile level $q\in(0,1)$, let $\widehat{\boldsymbol y}^{q}$ be its vector of predictions on a tuning sample, let $\operatorname{RMSE}$ be root mean squared error over that sample, and let $\operatorname{LateRMSE}$ and $\operatorname{EarlyRMSE}$ be the same error restricted to observed late and early deliveries. Given nonnegative weights $\beta$ on late errors and $\lambda$ on early errors, the selected quantile level $q^*$ is
\[
q^*\in\argmin_{q \in (0, 1)}
\operatorname{RMSE}(\widehat{\boldsymbol y}^q)
+\beta\operatorname{LateRMSE}(\widehat{\boldsymbol y}^q)
+\lambda\operatorname{EarlyRMSE}(\widehat{\boldsymbol y}^q).
\]
Changing these weights shifts the point forecast without refitting the predictive distribution. In implementation, the objective is evaluated over a grid of quantile levels on the validation data, and the selected $q^*$ is then applied to the test distributions.

\section{Contextual Fulfillment Optimization}
\label{sec:optimization-framework}

\begin{figure}[pos=!ht]
    \centering
    \includegraphics[width=\textwidth]{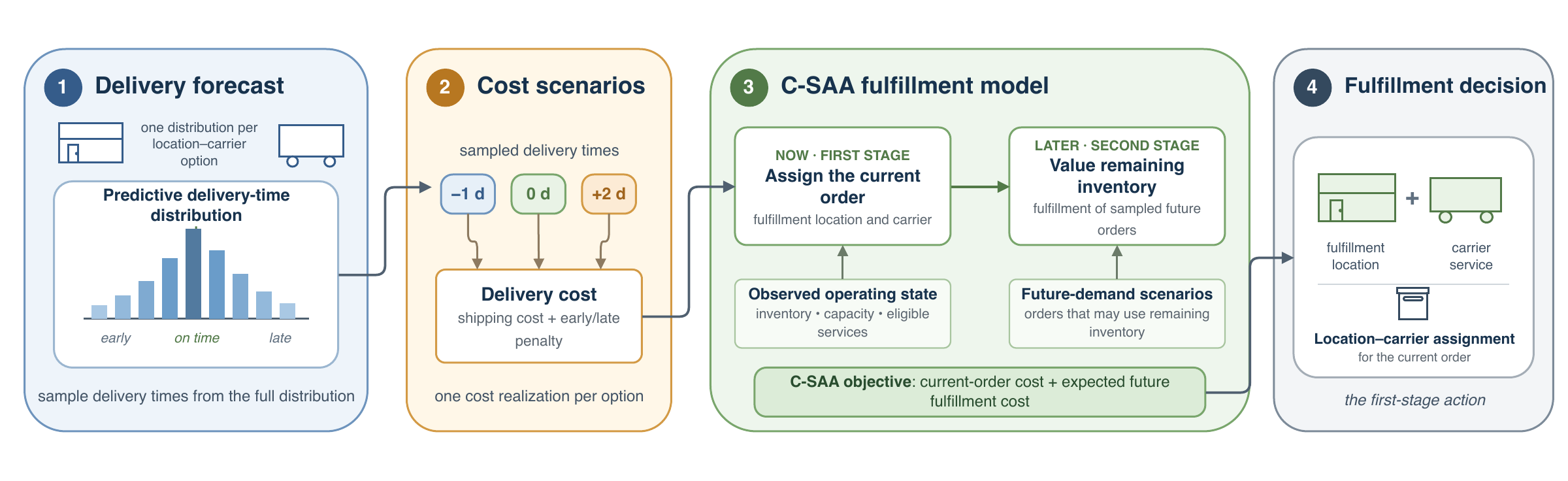}
    \caption{Prediction-to-optimization interface for fulfillment-location and carrier selection.}
    \label{fig:prediction-optimization}
\end{figure}

Fulfillment decisions balance shipping costs, consolidation savings, delivery penalties, and the cost of leaving demand unfulfilled, subject to inventory constraints. Because each decision consumes inventory, it also affects the resources available for later orders. Fulfillment options with the same point forecast can still differ in expected early- and late-delivery penalties. C-SAA therefore samples the forecast distributions rather than point forecasts to construct delivery-cost scenarios and combines these with future-demand scenarios when evaluating feasible decisions. Figure \ref{fig:prediction-optimization} summarizes how the scenarios enter the fulfillment decision.

\FloatBarrier

\subsection{Rolling C-SAA Policy}
\label{sec:csaa-model}

Let $\mathcal V_t(\boldsymbol s_t)$ be the minimum expected cost from epoch $t$ to the end of the horizon in problem \eqref{eq:sequential-problem}, given state $\boldsymbol s_t$ at epoch $t$, with $\mathcal V_{T+1}\equiv0$. The optimal costs satisfy the Bellman equation
\begin{equation}
\mathcal V_t(\boldsymbol s_t)=\min_{\boldsymbol z\in\mathcal Z(\boldsymbol s_t)}
\Bigl\{\E\bigl[c(\boldsymbol s_t,\boldsymbol z,\boldsymbol Y_t)\mid \boldsymbol s_t\bigr]
+\E\bigl[\mathcal V_{t+1}\bigl(f(\boldsymbol s_t,\boldsymbol z,O_{t+1})\bigr)\mid \boldsymbol s_t\bigr]\Bigr\}.
\label{eq:bellman}
\end{equation}
Here $\boldsymbol z$ is a candidate fulfillment decision for order $t$, and $\boldsymbol Y_t$ collects the delivery-time deviations of the eligible fulfillment options. The first term is the expected cost of the current order under $\boldsymbol z$, and the second is the expected cost of later orders given the state that $\boldsymbol z$ leaves. Because $\boldsymbol Y_t$ does not affect later states, the delivery-time distributions of the current order do not enter the second term. The second term depends on the inventory of every item at every fulfillment location and is intractable at the scale of the application, so the method follows the rolling two-stage approximation of \citet{ye2026proxy}.

At each epoch, contextual sample-average approximation (C-SAA) builds scenarios for two uncertain quantities: the delivery-time deviation of each eligible fulfillment option, drawn from the forecasts of Section \ref{sec:forecast-methodology}, and the future demand for the items in the current order, drawn from the probabilistic item-level demand forecast of \citet{ye2026proxy}. It then minimizes the average cost over these scenarios. The first stage chooses the fulfillment decision for the current order, and the second allocates the remaining inventory to later demand. Only the first-stage decision is implemented, and repeating the procedure for each new order defines the rolling C-SAA policy.

Let $\mathcal I_t$ contain the items in the order arriving at epoch $t$. For each item $i\in\mathcal I_t$, let $\mathcal A_i(\boldsymbol s_t)$ be its eligible fulfillment options, and let $\mathcal A(\boldsymbol s_t)=\cup_{i\in\mathcal I_t}\mathcal A_i(\boldsymbol s_t)$. The first-stage variable $z_{ia}$ is the number of units of item $i$ assigned to fulfillment option $a$, and $w_i$ is the number of its units left unfulfilled. The order-level decision $\boldsymbol z$ collects both, and the feasible set $\mathcal Z(\boldsymbol s_t)$ enforces order balance, inventory, and eligibility constraints. Let $u_a(\boldsymbol z)=\sum_{i:a\in\mathcal A_i(\boldsymbol s_t)}z_{ia}$ be the units assigned to fulfillment option $a$, $\operatorname{Disc}(\boldsymbol z)$ the consolidation savings generated by $\boldsymbol z$, and $\rho$ the lost-sales penalty.

Let $\kappa_{a,t}$ be the base shipping cost and $\gamma^-$ and $\gamma^+$ the early- and late-delivery penalty rates. The per-unit delivery cost at deviation $y$ is $c_{a,t}(y)=\kappa_{a,t}+\gamma^-\max\{-y,0\}+\gamma^+\max\{y,0\}$. The first-stage order cost $R_t(\boldsymbol z;\boldsymbol s_t)=\rho\sum_{i\in\mathcal{I}_t}w_i-\operatorname{Disc}(\boldsymbol z)$ charges unfulfilled units and credits consolidation savings, so the current-order cost of Section \ref{sec:problem-setting} is $c(\boldsymbol s_t,\boldsymbol z,\boldsymbol Y_t)=R_t(\boldsymbol z;\boldsymbol s_t)+\sum_{a\in\mathcal A(\boldsymbol s_t)}c_{a,t}(Y_{a,t})u_a(\boldsymbol z)$. Let $\boldsymbol D_t$ denote the total demand for the items in the current order from the orders that arrive later in the planning horizon. The recourse cost $Q_t(\boldsymbol z,\boldsymbol D_t;\boldsymbol s_t)$ is the minimum future shipping and shortage cost when that demand is served with the inventory that $\boldsymbol z$ leaves. The recourse allocation uses the constraints of \citet{ye2026proxy} and charges future shipments their base shipping cost. Neither $R_t$ nor $Q_t$ depends on the delivery-time outcomes.

Under the two-stage approximation, the expected cost of decision $\boldsymbol z$ at epoch $t$ is
\begin{equation}
\label{eq:population-objective}
J(\boldsymbol z;\boldsymbol s_t)=R_t(\boldsymbol z;\boldsymbol s_t)+\E\left[
  \sum_{a\in\mathcal A(\boldsymbol s_t)}c_{a,t}(Y_{a,t})u_a(\boldsymbol z)
  +Q_t\bigl(\boldsymbol z,\boldsymbol D_t;\boldsymbol s_t\bigr)\right],
\end{equation}
where the expectation is over each $Y_{a,t}$ given $\boldsymbol x_{a,t}$ and over future demand, which follows the demand forecast of \citet{ye2026proxy}.
Equation \eqref{eq:population-objective} is the expression minimized in the Bellman equation \eqref{eq:bellman} with the expected cost of later orders replaced by the expected recourse cost. The expected current-order cost is unchanged, since $R_t(\boldsymbol z;\boldsymbol s_t)+\sum_{a\in\mathcal A(\boldsymbol s_t)}\E\bigl[c_{a,t}(Y_{a,t})\mid\boldsymbol x_{a,t}\bigr]u_a(\boldsymbol z)=\E\bigl[c(\boldsymbol s_t,\boldsymbol z,\boldsymbol Y_t)\mid \boldsymbol s_t\bigr]$. The recourse cost simplifies the future in three ways: it allocates the remaining inventory in one step to the aggregate demand $\boldsymbol D_t$ after observing it, rather than order by order as orders arrive; it charges no delivery penalties on future shipments; and it ignores the effect of $\boldsymbol z$ on fulfillment-location workload. It retains the inventory opportunity cost of the current decision, because the inventory that $\boldsymbol z$ leaves at each fulfillment location determines whether later demand is shipped from more expensive locations or left unfulfilled. Minimizing $J(\cdot\,;\boldsymbol s_t)$ at each epoch and implementing the minimizer defines a policy for the sequential problem.

The delivery-time distributions are unknown, so C-SAA replaces each $F(\cdot\mid\boldsymbol x_{a,t})$ in Equation \eqref{eq:population-objective} by its forecast $\widehat F(\cdot\mid\boldsymbol x_{a,t})$. The result is the forecast objective $\widehat J(\boldsymbol z;\boldsymbol s_t)$, which C-SAA approximates by an average over sampled scenarios. Each scenario $\omega$ pairs a future-demand scenario $\boldsymbol D_t^{(\omega)}$ from the demand forecast with a delivery-time deviation $Y_{a,t}^{(\omega)}$ for each $a\in\mathcal A(\boldsymbol s_t)$, drawn from $\widehat F(\cdot\mid \boldsymbol x_{a,t})$. Sampling each fulfillment option from its own conditional distribution, rather than from a pooled historical sample, makes the approximation contextual: two fulfillment options for the same order can differ in both their expected delivery time and its spread. With $N$ sampled scenarios, indexed by $\omega=1,\ldots,N$, C-SAA solves
\begin{align}
\label{eq:forecast-saa}
\min_{\boldsymbol z\in\mathcal Z(\boldsymbol s_t)}\quad
    & R_t(\boldsymbol z;\boldsymbol s_t)
      +\frac{1}{N}\sum_{\omega=1}^{N}\left[
      \sum_{a\in\mathcal A(\boldsymbol s_t)}c_{a,t}(Y_{a,t}^{(\omega)})u_a(\boldsymbol z)
      +Q_t\bigl(\boldsymbol z,\boldsymbol D_t^{(\omega)};\boldsymbol s_t\bigr)
      \right].
\end{align}
The first-stage decision is shared across scenarios, and the recourse allocation varies with future demand. To reduce sensitivity to one scenario sample, the model solves $L$ independent instances of \eqref{eq:forecast-saa}, each with $N=N_1$ independently drawn scenarios \citep{kleywegt2002sample}. The candidates returned as optimal within the time limit are then evaluated on a common independent sample of $N_2$ scenarios, and the one with the lowest sample-average objective is implemented.
\FloatBarrier

\subsection{Theoretical Properties}
\label{sec:theory-guarantees}

In C-SAA, a forecast matters through the delivery-time scenarios drawn from it, and with enough scenarios its effect on the decision reduces to the expected delivery cost it implies at each fulfillment option. The validity guarantees of Section \ref{sec:conformal} concern the coverage of CPS prediction intervals and the calibration of binary IVAP probabilities. Both hold at the population level, or for MCPS within a Mondrian group, whereas a fulfillment decision depends on the expected costs of an individual order. A forecast that gives each fulfillment option the same delivery-time distribution for every order can satisfy these guarantees exactly, yet its decisions can exceed the optimal expected cost by an amount that grows with the late-delivery penalty (Appendix \ref{app:marginal-calibration}). The results below therefore bound the regret of the C-SAA decision for an individual order by the error in these expected costs, and then bound that error for CPS and for class-probability forecasts.

Fix the state $\boldsymbol s_t$ at which order $t$ arrives, with a finite, nonempty set of feasible decisions and eligible fulfillment options. Let $\boldsymbol z^*$ minimize the expected cost $J(\cdot\,;\boldsymbol s_t)$ of \eqref{eq:population-objective} over $\mathcal Z(\boldsymbol s_t)$, which is the best decision for the current order and inventory state under the true delivery-time distributions. The regret of a decision $\boldsymbol z$ is $J(\boldsymbol z;\boldsymbol s_t)-J(\boldsymbol z^*;\boldsymbol s_t)$, the amount by which its expected cost exceeds that of $\boldsymbol z^*$. Because $J$ takes future demand from the same demand forecast for every delivery-time forecast, the regret does not include demand-forecast error. The benchmark is optimal for one decision epoch within the two-stage approximation, not for the sequential problem. C-SAA selects $\widehat{\boldsymbol z}$ on the $N_2$ evaluation scenarios from a nonempty set $\mathcal C(\boldsymbol s_t)\subseteq\mathcal Z(\boldsymbol s_t)$ of candidates retained from its $L$ sample-average problems, and restricting the choice to these candidates adds the candidate-generation gap
\(
\Gamma_{\mathcal C}^J=\min_{\boldsymbol z\in\mathcal C(\boldsymbol s_t)}J(\boldsymbol z;\boldsymbol s_t)-J(\boldsymbol z^*;\boldsymbol s_t).
\)
The regret of the C-SAA decision therefore combines forecast error in the expected delivery costs, sampling error in the candidate evaluation, and $\Gamma_{\mathcal C}^J$.

\begin{samepage}
\begin{proposition}[Regret bound for any forecast distribution]
\label{thm:forecast-decision-loss}
Suppose the true and forecast objectives share the same first-stage and recourse costs, that each feasible decision assigns nonnegative units totalling at most $\bar u$, that $\delta\geq0$ bounds the per-unit error in expected delivery cost at every eligible fulfillment option, and that $\varepsilon\geq0$ bounds the gap between each retained candidate's evaluation-sample objective and its expected objective under the forecast distribution. Then
\[
J(\widehat{\boldsymbol z};\boldsymbol s_t)-J(\boldsymbol z^*;\boldsymbol s_t)
\leq 2\bar u\delta+2\varepsilon+\Gamma_{\mathcal C}^J.
\]
\end{proposition}
\end{samepage}

Because $\delta$ measures error in expected cost, a forecast with the wrong shape but accurate expected costs adds no forecast error. The factor of two arises because the same error can make the selected decision look cheaper than it is and the best retained candidate look more expensive, and $\varepsilon$ enters twice for the same reason. With $\varepsilon=0$ and $\Gamma_{\mathcal C}^J=0$, the bound has the same form as the plug-in regret bound of \citet[Theorem 4]{hu2022fast} at a single context.

A CPS forecast shifts the regression prediction by residuals from held-out calibration orders, so its accuracy depends on how well those residuals represent the residual at the current order. SCPS uses all of these residuals, and MCPS uses those from orders with similar predicted values. Let $b(a)$ be the calibration group used by fulfillment option $a$. The expected-cost error of the forecast has three sources. First, each group's empirical residual distribution differs from its population counterpart by sampling error. When the population residuals are conditionally iid within fixed calibration groups, the Dvoretzky--Kiefer--Wolfowitz inequality \citep{dvoretzky1956asymptotic,massart1990tight} bounds this error for all groups at once by a term that decreases with the group size. Second, the calibration residuals come from the fitted model rather than the population regression function. This adds $\zeta^{\rm cal}_{b(a)}$, the largest absolute regression error on the group's calibration observations, and $\zeta^{\rm te}_a$, the absolute regression error at the current input, each scaled by the Lipschitz constant $M_{b(a)}$ of the group's population residual CDF. Third, the residual distribution at the current order may differ from that of its calibration group, and $\xi_{a,\boldsymbol s_t}$ denotes the largest absolute difference between their CDFs. These three terms bound the CDF error of the regression response, and the total variation $V_{a,t}$ of the delivery cost, as a function of that response, converts them into expected-cost error. Finally, a Hoeffding bound $e_{N_2}$ (Equation \eqref{eq:common-sample-error} in Appendix \ref{app:theory-setup}) replaces $\varepsilon$. It holds for all retained candidates with probability at least $1-\alpha_s$, for a chosen $\alpha_s\in(0,1)$, and decreases as $N_2^{-1/2}$. Appendix \ref{app:theory-setup} states these quantities and their assumptions in full.

\begin{proposition}[Regret bound for the CPS residual sampler]
\label{thm:cps-csaa}
Under the conditions of Appendix \ref{app:theory-setup}, let $B$ be the number of calibration groups across all fitted forecasting models, one per model for SCPS and one per model and prediction range for MCPS, let group $b(a)$ contain $N_c^{b(a)}$ calibration observations, and let $\alpha_c\in(0,1)$. Define
\[
e_{b(a)}=\sqrt{\frac{\log(2B/\alpha_c)}{2N_c^{b(a)}}},
\qquad
\delta^{\rm CPS}_{a,t}=V_{a,t}\bigl[e_{b(a)}+M_{b(a)}\bigl(\zeta^{\rm cal}_{b(a)}+\zeta^{\rm te}_a\bigr)+\xi_{a,\boldsymbol s_t}\bigr].
\]
Then, with probability at least $1-\alpha_c-\alpha_s$ over the calibration and scenario samples,
\begin{equation}
J(\widehat{\boldsymbol z};\boldsymbol s_t)-J(\boldsymbol z^*;\boldsymbol s_t)
\leq 2\bar u\max_{a\in\mathcal A(\boldsymbol s_t)}\delta^{\rm CPS}_{a,t}
+2e_{N_2}+\Gamma_{\mathcal C}^J.
\label{eq:cps-costgap}
\end{equation}
\end{proposition}
Proposition \ref{thm:cps-csaa} shows the trade-off between SCPS and MCPS. Pooling all residuals reduces the sampling term $e_{b(a)}$, whereas grouping orders with similar predicted values can reduce the mismatch $\xi_{a,\boldsymbol s_t}$. The proposition takes the Mondrian groups as fixed, which excludes bin boundaries chosen from the calibration data (Appendix \ref{app:theory-proof}).

The bound does not use the rank-based validity of Section \ref{sec:cps}. The conformal construction enters instead through its held-out calibration set, which keeps the model fixed given the proper training set, so $\zeta^{\rm cal}_{b(a)}$ measures out-of-sample regression error. Residual-based contextual SAA builds scenarios in the same way from in-sample or leave-one-out residuals, and its analysis by \citet{kannan2022data} assumes that regression errors differ across contexts only in scale and that cost is Lipschitz in the outcome. Proposition \ref{thm:cps-csaa}, by contrast, lets the residual distribution change in any way across orders at the price of the mismatch term $\xi_{a,\boldsymbol s_t}$, and allows discontinuous costs, which arise when a continuous delivery-time response is converted to whole days.

A class-probability forecast such as MC-CLF assigns probabilities to a finite set of delivery outcomes, and the delivery cost takes one value in each outcome class. Its expected-cost error therefore depends on how far the forecast class probabilities are from the true ones at the current order. For fulfillment option $a$, let $\boldsymbol\theta_a=(\theta_{ak})_{k=1}^K$ and $\boldsymbol p_a=(p_{ak})_{k=1}^K$ be the true and forecast class-probability vectors over the $K$ outcome classes, and let $\Delta_a$ be the range of delivery costs across those classes. The excess conditional Brier score and log loss measure the distance between the two vectors.

\begin{proposition}[Regret bound for class-probability forecasts]
\label{thm:va-csaa}
Under the finite-support and evaluation conditions of Appendix \ref{app:theory-setup}, define
\[
r^{\rm B}_a=\|\boldsymbol p_a-\boldsymbol\theta_a\|_2^2,
\qquad
r^{\log}_a=\sum_{k=1}^K \theta_{ak}\log\!\left(\frac{\theta_{ak}}{p_{ak}}\right),
\qquad
\delta_{a}=\frac{\Delta_a}{2}\min\!\left\{\sqrt{K r^{\rm B}_a},\sqrt{2r^{\log}_a}\right\},
\]
with terms at $\theta_{ak}=0$ contributing zero and $r^{\log}_a$ infinite when $p_{ak}=0<\theta_{ak}$. Then, for any normalized class-probability forecast, with probability at least $1-\alpha_s$ over the scenario sample,
\begin{equation}
J(\widehat{\boldsymbol z};\boldsymbol s_t)-J(\boldsymbol z^*;\boldsymbol s_t)
\leq 2\bar u\max_{a\in\mathcal A(\boldsymbol s_t)}\delta_{a}
+2e_{N_2}+\Gamma_{\mathcal C}^J.
\label{eq:va-costgap}
\end{equation}
\end{proposition}
The Brier bound is a single-order counterpart of the squared-loss risk bound of \citet{honguyen2022risk}. The cost range $\Delta_a$ sets how much a probability error costs, so accuracy matters most at fulfillment options whose delivery outcomes differ widely in cost. Venn--Abers calibration does not by itself control $r^{\rm B}_a$. Averaged over orders, the excess Brier score splits into a calibration loss, which Venn--Abers calibration targets, and a grouping loss from differences among orders that receive the same forecast \citep{kull2015novel}. Recalibration cannot reduce the grouping loss because it assigns the same probabilities to orders with the same forecast. For binary decisions, this loss still drives regret after recalibration \citep{perezlebel2025decision}. For a Venn--Abers forecast, the width of its probability intervals gives an alternative bound that can be computed from the forecast (Lemma \ref{lem:va-width} in Appendix \ref{app:theory-proof}), but it assumes that the intervals contain the true class probabilities at the current order, which Venn--Abers validity does not imply. Appendix \ref{app:theory-proof} gives the proofs.

\section{Forecasting Study}
\label{sec:experiments}
This study evaluates how well each method represents delivery-time uncertainty on proprietary order records. 

\subsection{Data and Experimental Setup}
\label{sec:industrial-data}
\label{sec:forecast-setup}

The forecast-level study uses an anonymized e-commerce transactional dataset provided by an industrial partner. It contains more than five million ship-to-home orders collected between July 2023 and November 2024 and is commercially sensitive, so it cannot be released. The analysis focuses on a ground carrier service selected with the industrial partner and retains delivery-time deviations from $-10$ to $10$ days, which limits the influence of data errors and atypical cases. Within this range, the deviation is centered near zero but has substantial early and late tails. Both tails matter operationally, since late delivery incurs service penalties, whereas early arrival can disrupt customer and downstream planning \citep{salari2022real,cui2024sooner}.

\begin{figure}[pos=!ht]
\centering
    \includegraphics[width=\textwidth]{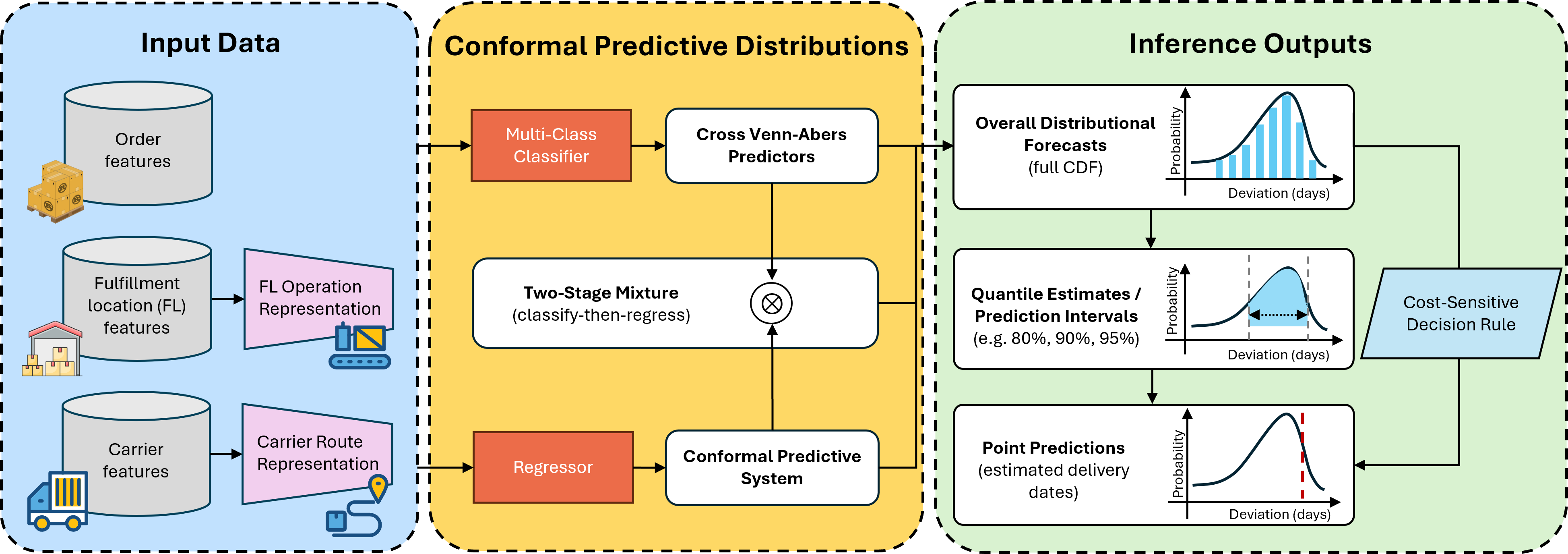}
    \caption{Forecasting framework used in the proprietary-data study.}
    \label{fig:framework}
\end{figure}

For the proprietary dataset, the forecasting models use order, fulfillment-location, and carrier features. Order features describe the shipment and customer, fulfillment-location features summarize workload, congestion, and processing time, and carrier features describe route geography, transit stages, tracking activity, recent performance, and delay exceptions. Figure \ref{fig:framework} summarizes how these inputs enter the models and how the resulting distributions are used. Relative to \citet{ye2025msom}, this dataset adds variables such as customer membership tier and timestamps describing fulfillment-location and carrier-network operations. Appendix Table \ref{tab:all_features} lists the complete feature set, and Appendix Figure \ref{fig:feat-imp} reports XGBoost feature importance, measured by each feature's contribution to reducing the tree-training loss. The two most important features, the number of days from order to promised delivery and the shipping-offer method at order time, both describe the delivery promise from which the deviation is measured.

\begin{figure}[pos=!ht]
    \centering
    \includegraphics[width=0.6\textwidth]{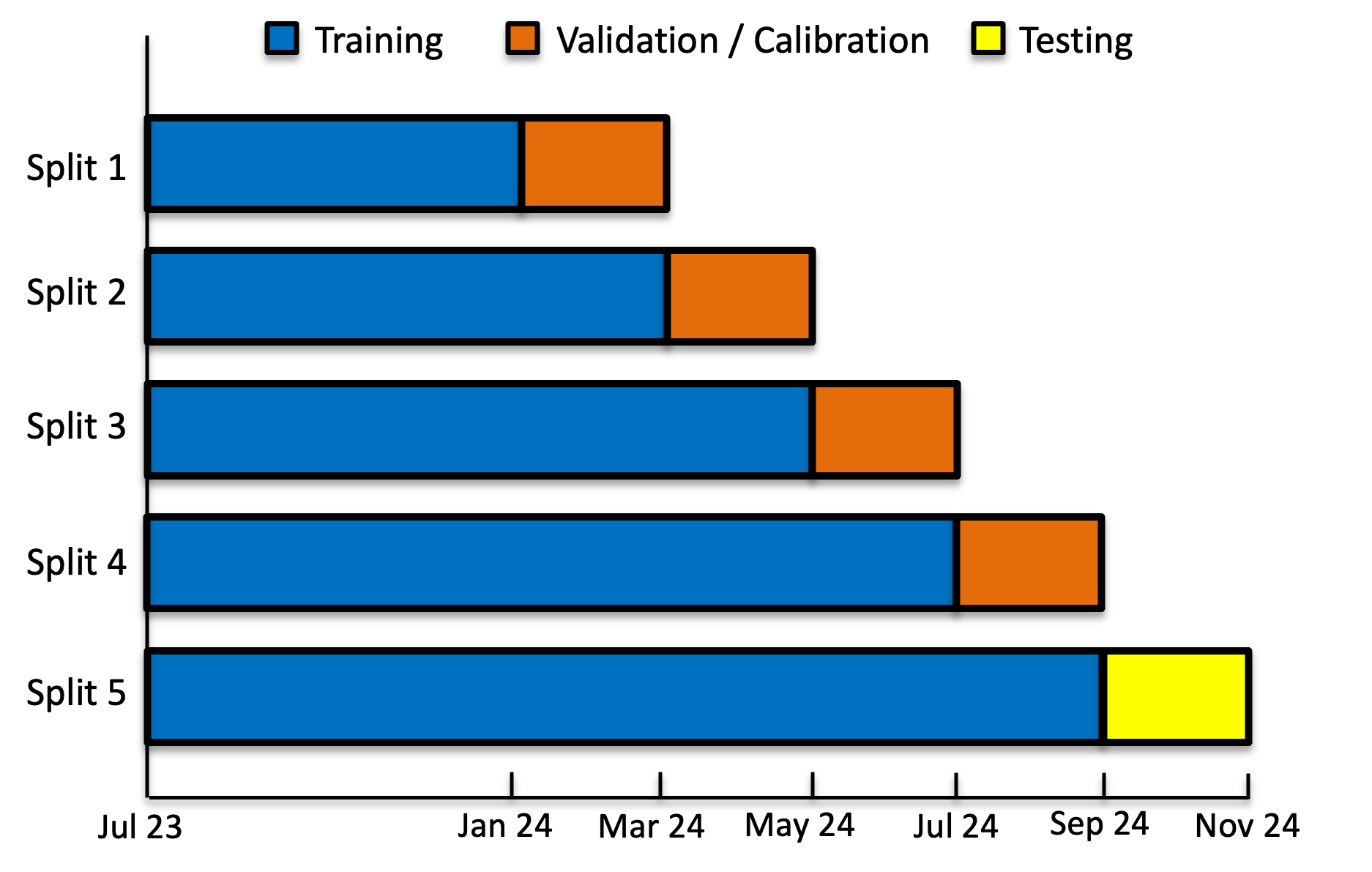}
    \caption{Data Splitting for the Experimental Results.}
    \label{fig:data-split}
\end{figure}

\paragraph{Temporal data split.}
The forecast study follows the chronological protocol of \citet{ye2025iccl}, which keeps later observations out of model fitting and calibration and evaluates every method on the same final test period. Figure \ref{fig:data-split} shows the five chronological data splits. In Splits 1--4, the training period expands forward in time and the next period supplies validation and calibration observations for model selection. Split 5 holds out the final test period. CVAP uses calibration observations from all four validation folds. SCPS and MCPS use the calibration observations in the final validation fold. The test set consists of orders shipped by the carrier service that the existing system evaluated when the fulfillment location was selected, so the machine-learning forecasts and the existing system's estimate refer to the same service.

\FloatBarrier

\paragraph{Implementation.}
XGBoost \citep{chen2016xgboost} is the base regressor for CPS, and CatBoost \citep{prokhorenkova2018catboost} is the base MC-CLF model. The weighted-error point forecasts use $\beta=0.5$ and $\lambda=0$, placing additional weight on errors for late deliveries. The continuous ranked probability score (CRPS) evaluates the complete predictive distribution. Pinball loss (PL) evaluates quantile accuracy, and mean quantile coverage error (MQCE) measures the deviation between empirical and nominal quantile coverage. Point forecasts are evaluated by exact-day accuracy, RMSE, and mean absolute error (MAE). Appendix \ref{app:computational-reproducibility} reports the software, tuning, and MCPS partition details, and Appendix Table \ref{tab:evaluation_metrics} gives the metric formulas.

\subsection{Forecasting Results}
\label{sec:results}
\label{sec:forecast-results}

The evaluation covers the full predictive distributions, which the fulfillment model uses, and the single delivery dates derived from them.

\subsubsection{Distributional Forecasting Performance}
Table \ref{tab:distribution_metrics} reports PL, MQCE, CRPS, and the coverage and mean width of the 80\%, 90\%, and 95\% prediction intervals. Besides the conformal forecasts, it includes the uncalibrated classifier (Base MC-CLF) and two benchmark distributional forecasts from \citet{ye2025msom}, the Quantile Random Forest (QRF) and isotonic-regression-calibrated MC-CLF (IR MC-CLF).

\begin{table}[!ht]
\centering
\scriptsize
\caption{Average Pinball Loss, MQCE, CRPS, Coverage, and Interval Sizes on the Test Set.}
\label{tab:distribution_metrics}
\begin{adjustbox}{width=0.95\textwidth}
\begin{tabular}{l|ccc|ccc|ccc}
\toprule
\textbf{Methods}
& \textbf{PL $\downarrow$} & \textbf{MQCE $\downarrow$} & \textbf{CRPS $\downarrow$}
& \multicolumn{3}{c|}{\textbf{Coverage}} & \multicolumn{3}{c}{\textbf{Mean Interval Size $\downarrow$}} \\
&  &  &
& 80\% & 90\% & 95\% & 80\% & 90\% & 95\% \\
\midrule
Base MC-CLF & 0.335 & 0.232 & 0.498 & 0.427 & 0.600 & 0.768 & 2.043 & 2.749 & 3.517 \\
IR MC-CLF   & 0.298 & 0.240 & \textbf{0.413} & 0.619 & 0.776 & 0.901 & 2.232 & 2.944 & 3.659 \\
CVAP MC-CLF & 0.297 & 0.232 & \textbf{0.419} & 0.658 & 0.835 & 0.935 & 2.312 & 3.049 & 3.780 \\
\midrule
QRF         & 0.273 & 0.313 & 0.502 & 0.876 & 0.923 & 0.965 & \textbf{1.844} & \textbf{2.587} & \textbf{3.445} \\
\midrule
SCPS        & 0.261 & \textbf{0.035} & 0.482 & 0.817 & 0.907 & 0.953 & 2.133 & 3.374 & 4.582 \\
MCPS        & \textbf{0.244} & \textbf{0.041} & 0.451 & 0.829 & 0.915 & 0.953 & 2.097 & 3.002 & 4.200 \\
\midrule
2-stage SCPS & \textbf{0.238} & 0.217 & \textbf{0.440} & 0.795 & 0.879 & 0.935 & \textbf{1.439} & \textbf{2.264} & \textbf{3.172} \\
2-stage MCPS & \textbf{0.242} & \textbf{0.213} & 0.449 & 0.828 & 0.886 & 0.925 & \textbf{1.664} & \textbf{2.391} & \textbf{3.135} \\
\bottomrule
\multicolumn{10}{l}{\scriptsize\textit{Note: $\downarrow$ indicates lower values are better. Best three values are \textbf{bolded}.}}\\
\end{tabular}
\end{adjustbox}
\end{table}

Both isotonic regression and CVAP improve the uncalibrated MC-CLF distribution and achieve similar overall forecast accuracy. SCPS and MCPS provide the most reliable quantiles: both attain near-nominal interval coverage and improve pinball loss relative to QRF. MCPS also narrows the intervals and lowers CRPS relative to SCPS with similar calibration, so grouping residuals by the regression prediction sharpens the distribution without losing coverage.

\begin{figure}[pos=!ht]
    \centering
    \includegraphics[width=\textwidth]{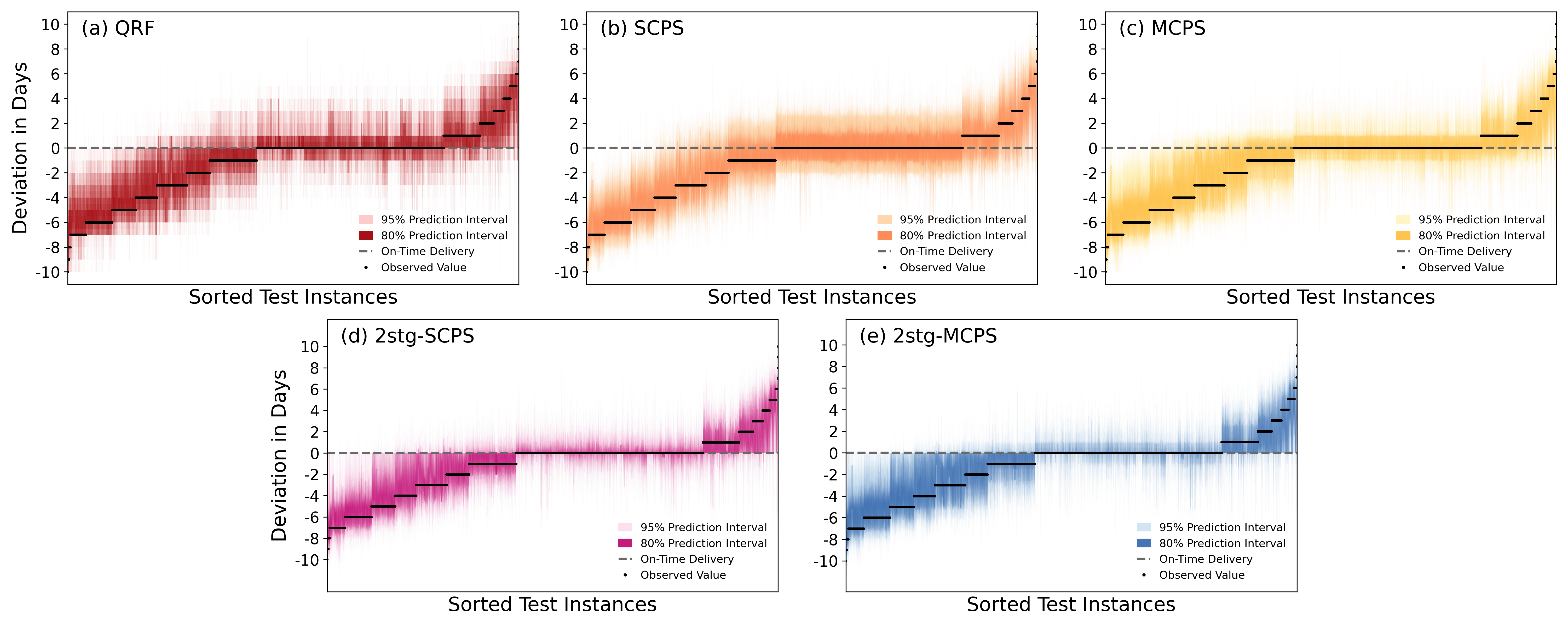}
    \caption{Prediction intervals on sorted test instances. The black marks indicate realized deviations. The label 2stg denotes the 2-stage SCPS and MCPS forecasts.}
    \label{fig:intervals}
\end{figure}

The 2-stage forecasts give the smallest mean interval widths in Table \ref{tab:distribution_metrics}. Their status probabilities and component-specific intervals distinguish the likelihood of a delay from its likely magnitude \citep{ye2025iccl}, and Figure \ref{fig:intervals} shows the resulting intervals for individual test orders. The narrower intervals, however, come with less consistent aggregate coverage and quantile calibration than those of SCPS and MCPS.

\FloatBarrier

\subsubsection{Point Prediction Performance}

The predictive distributions can also be reduced to a single delivery date and compared with the current system, whose estimate comes from transit tables and cutoff times. MC-CLF uses its most probable deviation class, and the base regression uses its rounded point estimate. For CPS, the reported point estimates are the median and the quantile selected by the weighted-error criterion. The same criterion, which emphasizes late errors, is also applied to CVAP MC-CLF and the 2-stage distributions.

\begin{table}[!ht]
\centering
\scriptsize
\centering
\caption{Point Prediction Performance on the Proprietary Test Set.}
\label{tab:point_metrics2}
\begin{adjustbox}{width=\textwidth}
\begin{tabular}{lcccccc}
\toprule
\textbf{Forecast / point summary} & \shortstack{\textbf{Overall}\\\textbf{accuracy (\%)} $\uparrow$}
& \shortstack{\textbf{Overall}\\\textbf{RMSE} $\downarrow$} & \shortstack{\textbf{Overall}\\\textbf{MAE} $\downarrow$}
& \shortstack{\textbf{Late-delivery}\\\textbf{detection (\%)} $\uparrow$}
& \shortstack{\textbf{Late}\\\textbf{RMSE} $\downarrow$} & \shortstack{\textbf{Late}\\\textbf{MAE} $\downarrow$}\\
\midrule
Current system & 49.53 & 1.66 & 1.01 & 3.82 & 2.99 & 2.49 \\
\midrule
Base MC-CLF & 58.50 & 1.30 & 0.71 & 48.67 & 1.84 & 1.31 \\
IR MC-CLF & \textbf{64.41} & 1.14 & 0.58 & 63.60 & 1.71 & 1.09\\
CVAP MC-CLF & \textbf{64.14} & 1.16 & 0.59 & 62.53 & 1.72 & 1.11\\
CVAP MC-CLF (weighted-error) & 52.47 & 1.30 & 0.77 & 73.60 & 1.62 & 1.09\\
\midrule
Base regression & 58.79 & \textbf{1.01} & \textbf{0.56} & \textbf{75.17} & \textbf{1.52} & \textbf{0.97} \\
SCPS (median) & 60.23 & \textbf{1.00} & \textbf{0.54} & 73.37 & 1.54 & 0.98\\
MCPS (median) & \textbf{61.47} & \textbf{1.00} & \textbf{0.53} & 63.21 & 1.60 & 1.07\\
SCPS (weighted-error) & 54.32 & 1.05 & 0.62 & \textbf{79.37} & \textbf{1.48} & \textbf{0.95}\\
MCPS (weighted-error) & 50.63 & 1.12 & 0.68 & \textbf{76.66} & \textbf{1.47} & \textbf{0.93}\\
\midrule
2-stage SCPS (weighted-error) & 58.83 & 1.08 & 0.60 & 66.71 & 1.56 & 1.04\\
2-stage MCPS (weighted-error) & 58.18 & 1.08 & 0.60 & 59.99 & 1.59 & 1.08\\
\bottomrule
\multicolumn{7}{l}{\footnotesize\textit{Note: Arrows indicate whether higher ($\uparrow$) or lower ($\downarrow$) values are better.}}\\
\multicolumn{7}{l}{\footnotesize\textit{Bold numbers indicate the top three values for each metric.}}
\end{tabular}
\end{adjustbox}
\end{table}

Table \ref{tab:point_metrics2} shows that every data-driven model improves on the current system across the reported metrics. The largest practical difference is in late-delivery detection, where the rule-based estimate misses most late deliveries and the machine-learning models identify a substantially larger share. Calibrated MC-CLF performs best on exact-day accuracy, while the CPS medians improve the regression model's overall error measures. The MCPS median gives the strongest central prediction among the regression-based summaries, but its late-delivery detection is weaker than that of the SCPS median.

\begin{figure}[pos=!ht]
    \centering
    \begin{minipage}[t]{0.45\textwidth}
        \centering
        \includegraphics[width=\textwidth]{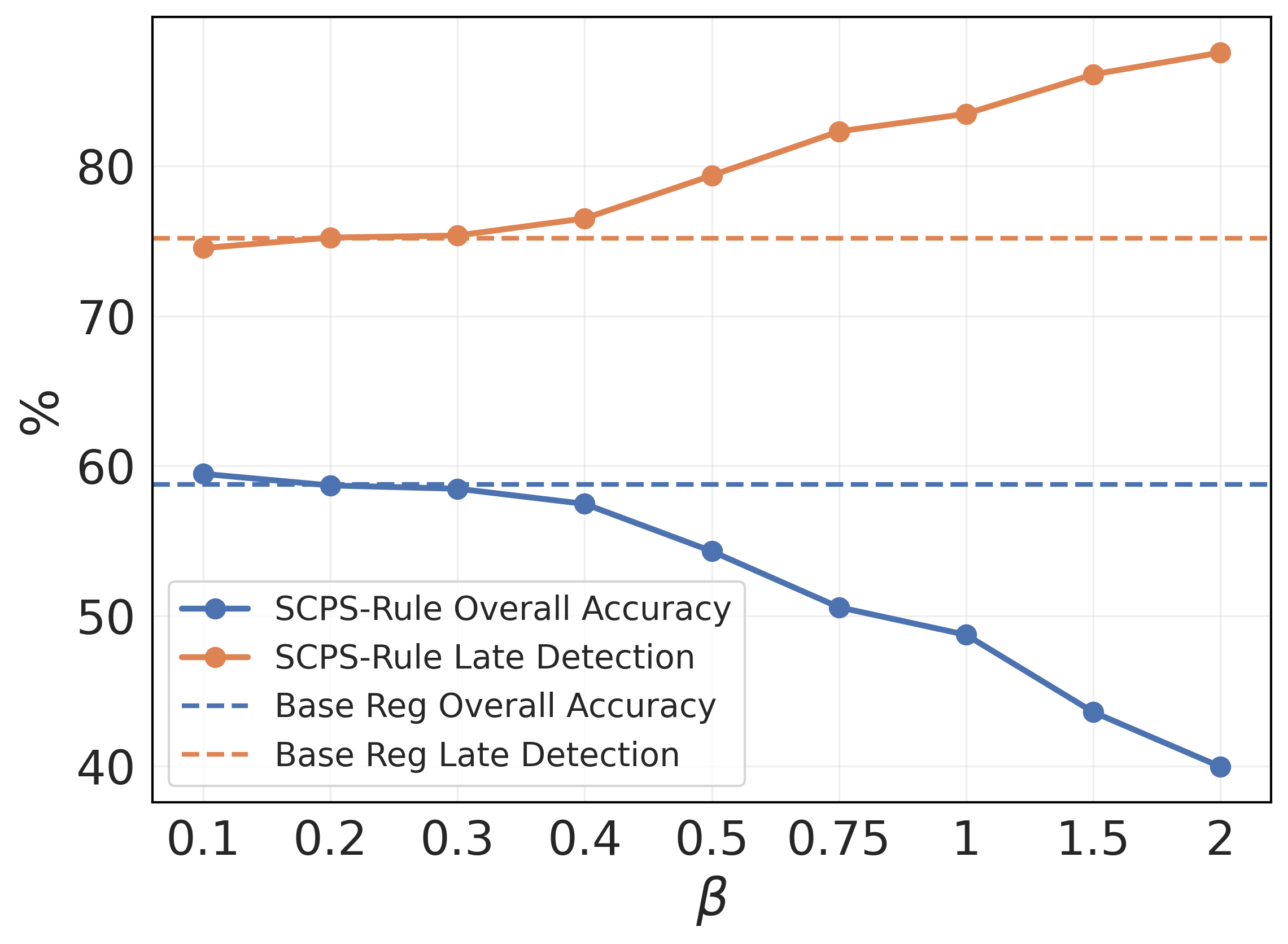}
    \end{minipage}
    \begin{minipage}[t]{0.45\textwidth}
        \centering
        \includegraphics[width=\textwidth]{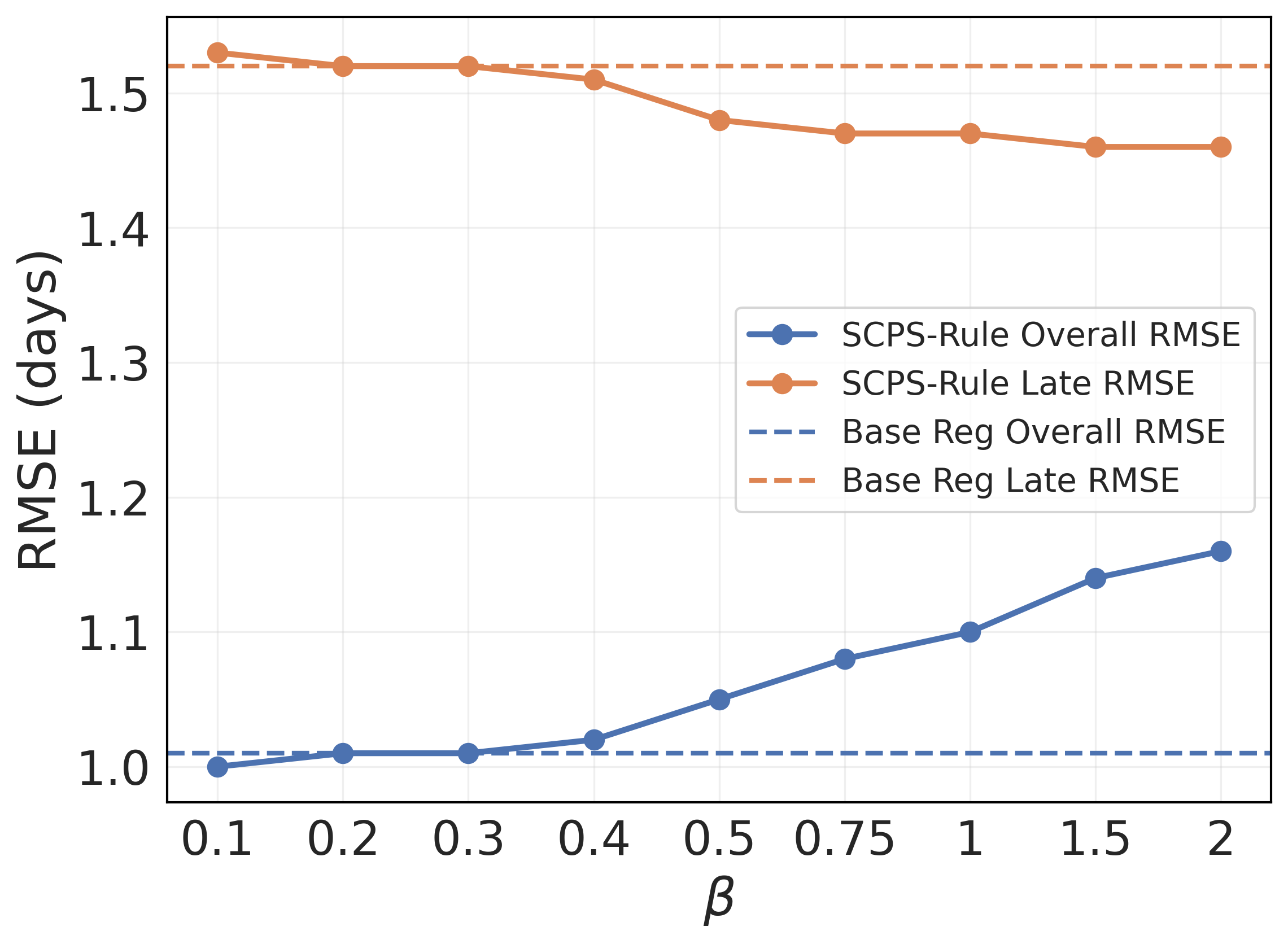}
    \end{minipage}
    \caption{Exact-day accuracy and RMSE of the SCPS point estimate versus the late-error weight $\beta$ in the point-prediction objective.}
    \label{fig:sensitivity_beta}
\end{figure}

The weighted-error estimates shift the trade-off toward late-delivery detection at the expense of exact-day accuracy. Raising the late-error weight $\beta$ continues that shift. In Figure \ref{fig:sensitivity_beta}, late-delivery detection and late RMSE improve as $\beta$ grows, until sufficiently large values reduce overall accuracy and raise overall RMSE above that of the base regression. A retailer that prioritizes identifying likely delays can choose a higher quantile, while one that prioritizes the most accurate single date can retain a central summary. The predictive distribution supports either use without retraining the model.

\FloatBarrier

\subsubsection{Operational Implications}
The forecasts support different operational uses. Because SCPS and MCPS intervals attain close to their nominal coverage, they can set customer delivery windows or the time windows used in last-mile route planning. The 2-stage forecasts separate the chance of a delay from its likely size, which suits simple triage rules, such as flagging orders whose probability of lateness exceeds a threshold for customer-service follow-up. Systems that accept a single date can use the MCPS median, the most accurate central estimate among the regression-based summaries. When missed delays are costly, they can use a weighted-error quantile instead, for example to send delay notices early or to schedule additional staff for shipments expected to arrive late.

\section{Decision Study}
\label{sec:decision-study}

The decision study asks whether differences among the forecast distributions carry through to fulfillment decisions. The proprietary records cannot answer this question, because each order records the delivery outcome only for the fulfillment option actually used, so the cost of any other decision is never observed. Evaluating fulfillment policies instead requires a simulator of orders, inventory, and delivery outcomes at every fulfillment option. The study therefore uses the simulator of \citet{ye2026proxy}, built from the public JD.com dataset of \citet{shen2024jd}, to compare the realized cost and late-delivery rate of the resulting C-SAA policies.

\subsection{Data and Experimental Setup}
\label{sec:jd-simulator}
\label{sec:decision-setup}

The simulator processes the same 38,416 peak-hour orders from six test dates in each of 50 replications. Within each replication, all policies share a new set of delivery outcomes. These outcomes are generated separately from the optimization scenarios by a model estimated on the test-period records. The forecasting models supplied to C-SAA are trained on JD.com records that precede the test dates.

C-SAA is run with five forecasts: the original nonconformal forecast, SCPS, MCPS, and MC-CLF calibrated by either the inductive (IVAP) or the cross (CVAP) Venn--Abers predictor. The nonconformal forecast is the original delivery-time model of \citet{ye2026proxy}, which supplies conditional delivery-time distributions without a conformal calibration step. SCPS and MCPS wrap a regression model with the same architecture as the nonconformal forecast, retrained on the proper-training period so that the calibration week is held out. MC-CLF is a CatBoost classifier over integer-day durations. Because these forecasts predict delivery duration rather than deviation, C-SAA converts each sampled duration to a deviation from the promised date before applying delivery penalties. Two benchmarks are forecast-free. Greedy assigns available inventory to fulfillment options in order of base shipping cost. Empirical SAA solves the stochastic program with delivery and demand scenarios sampled from historical observations. Appendix \ref{app:computational-reproducibility} reports forecast fitting, service coverage, cost coefficients, and the remaining simulation details.

The primary outcomes are total cost and the percentage of test orders delivered late. Total cost is the realized cost of all test orders over the six dates, including shipping costs, delivery penalties, and lost sales.

\subsection{Fulfillment Results}
\label{sec:decision-results}

Table \ref{tab:optimization_results} reports the mean total cost and the percentage of test orders delivered late across the six dates and replications. The confidence intervals quantify Monte Carlo variation conditional on this fixed six-date order trace and cost specification.

\begin{table}[!ht]
\centering
\footnotesize
\centering
\caption{Mean Total Cost and Late-Delivery Rate for C-SAA with Each Forecasting Method and for Forecast-Free Benchmarks.}
\label{tab:optimization_results}
\setlength{\tabcolsep}{6pt}
\begin{tabular*}{0.97\textwidth}{@{\extracolsep{\fill}}lrr@{}}
\toprule
\textbf{Policy / forecast} & \textbf{Mean total cost (\$)} $\downarrow$ & \textbf{Late-delivery rate (\%)} $\downarrow$ \\
\midrule
C-SAA + MC-CLF (CVAP) & 802,106 $\pm$ 3,246 & 3.93 $\pm$ 0.07 \\
C-SAA + MCPS & 805,419 $\pm$ 3,208 & 3.96 $\pm$ 0.07 \\
C-SAA + SCPS & 807,636 $\pm$ 3,563 & 4.11 $\pm$ 0.07 \\
C-SAA + MC-CLF (IVAP) & 817,676 $\pm$ 3,375 & 4.55 $\pm$ 0.07 \\
C-SAA + original nonconformal forecast & 863,422 $\pm$ 3,555 & 5.68 $\pm$ 0.07 \\
Greedy & 928,970 $\pm$ 4,066 & 9.80 $\pm$ 0.08 \\
Empirical SAA & 1,007,185 $\pm$ 3,551 & 11.28 $\pm$ 0.08 \\
\bottomrule
\end{tabular*}
\par\vspace{2pt}
\raggedright\scriptsize\textit{Note: Total cost is the six-date sum in each replication. Late-delivery rate is the percentage of test orders delivered late across the same six dates and 50 replications. Estimates are reported with 95\% confidence-interval half-widths.}
\end{table}

\begin{figure*}[pos=!ht]
    \centering
    \includegraphics[width=\textwidth]{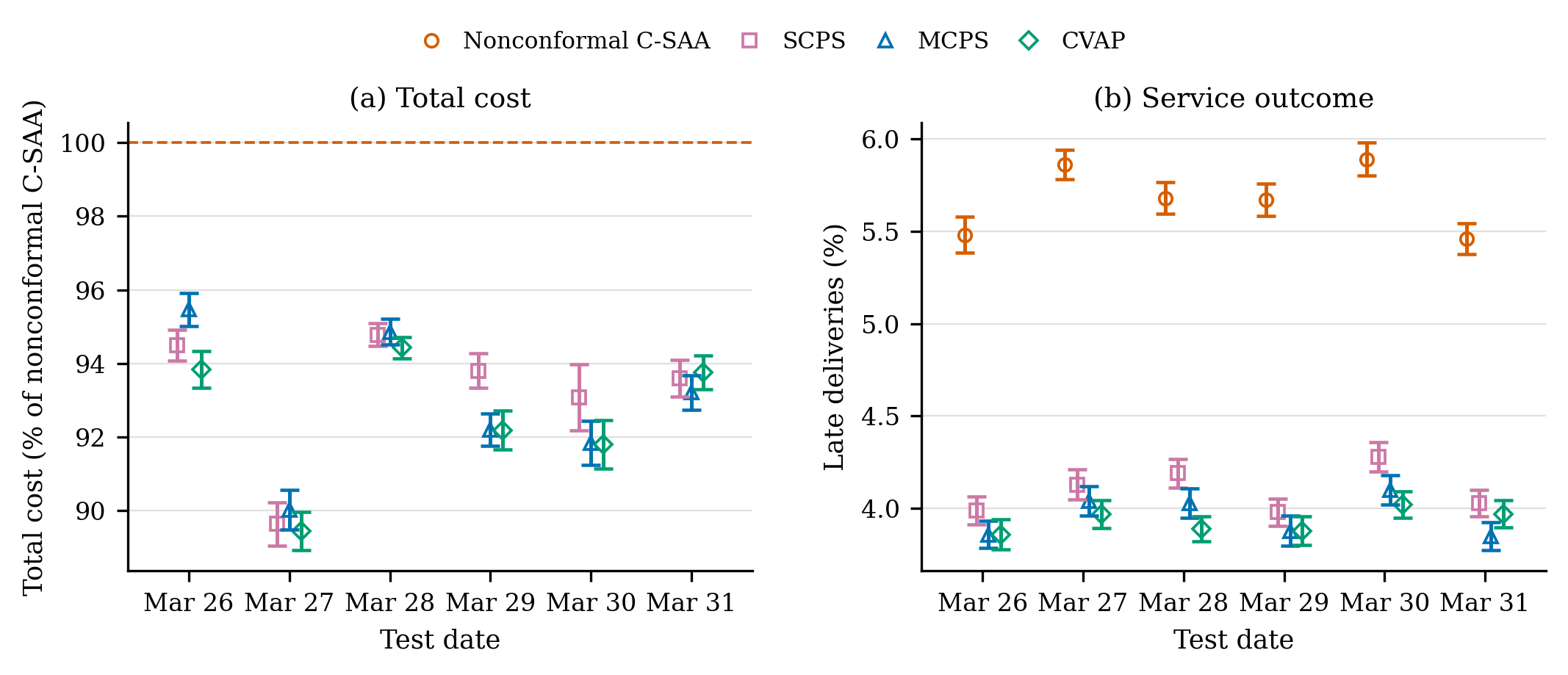}
    \caption{Daily performance of C-SAA with four forecasting methods. Panel (a) reports the total cost on each test date as a percentage of its mean under the original nonconformal forecast, shown by the dashed line, and panel (b) reports late-delivery rates. Markers are means over 50 replications, and error bars are 95\% confidence intervals for the means.}
    \label{fig:decision-across-dates}
\end{figure*}
\FloatBarrier

Replacing the forecast distribution alone changes the fulfillment decisions and their cost. Relative to C-SAA with the original nonconformal forecast, CVAP lowers the mean total cost by 7.1\% and MCPS by 6.7\%, while the optimization model and its solution procedure stay the same (Table \ref{tab:optimization_results}). The improvement holds on each of the six test dates, where CVAP, MCPS, and SCPS reduce the total cost by 4.5\% to 10.5\% and deliver fewer orders late (Figure \ref{fig:decision-across-dates}). Total cost under CVAP, MCPS, and SCPS also stays below that under the nonconformal forecast at every reported percentile (Appendix Table \ref{tab:optimization_quantiles}).

CVAP, MCPS, and SCPS perform similarly. MCPS and SCPS have overlapping confidence intervals, and CVAP differs from MCPS by less than 0.5\%. IVAP, which relies on a single calibration split, trails CVAP by a clear margin, and part of this gap may reflect the larger training set of CVAP, whose folds span the whole pre-test period. All machine-learning forecasts outperform the forecast-free benchmarks, consistent with the value of conditioning delivery-time uncertainty on the current order and fulfillment option.

The savings come from fewer late deliveries. Late unit-days, the number of units multiplied by their days of delay, fall from 0.137 per order under the nonconformal forecast to 0.101 under MCPS and 0.100 under CVAP. At the late-delivery penalty of \$40 per unit-day, this reduction accounts for 93\% of the cost saving under MCPS and CVAP and 83\% under SCPS, and shipping and early-delivery costs account for the rest. Lost sales are identical under these four forecasts, at 0.062 units per order. The conformal forecasts therefore change which fulfillment options receive the units, while leaving lost sales unchanged.

For a retailer, these results suggest that conformal calibration can be added to an existing delivery-time forecast without changing the optimization model. The choice between the two leading forecasts depends on what else the retailer needs. MCPS suits a retailer that also needs calibrated prediction intervals, for example to quote delivery dates, whereas CVAP suits one that already forecasts delivery-day classes. The cost differences rank the policies under the simulator and do not estimate savings in live operations, a limitation discussed in Section \ref{sec:conclusions}.

\section{Conclusions}
\label{sec:conclusions}

This paper evaluates conformal predictive distributions of delivery time both as forecasts and as inputs to contextual fulfillment optimization. The forecasting study shows how these distributions represent uncertainty on proprietary order records. The decision study samples delivery scenarios from CPS and CVAP distributions for a rolling C-SAA policy and examines how the distributions affect fulfillment choices in the JD.com simulator.

The two studies assess the forecasts in different ways. On the proprietary records, machine-learning forecasts improve delivery-date accuracy and late-delivery detection, while CPS provides well-calibrated prediction intervals. The 2-stage forecasts give narrower intervals and a delivery-status decomposition, with less consistent coverage. In the decision study, CVAP and MCPS give the lowest mean total cost among the machine-learning forecasts. Their advantage over the original nonconformal forecast shows that a conformal delivery-time distribution can improve decisions without any change to the optimization model or its solution procedure.

The study has limitations that point to future work. The decision study compares policies in a simulator, so its cost differences show how the policies rank rather than what they would save in practice. Measuring the actual effect requires a live deployment, such as an A/B test against the current policy, which would also show whether the advantage persists as conditions change. The theory also shows that the cost of a fulfillment decision depends on the accuracy of expected delivery costs for the current order, which marginal conformal guarantees do not control. Further work could calibrate forecasts for these expected costs, and extend the regret bounds of Section \ref{sec:theory-guarantees} from a single fulfillment decision to the sequence of decisions over the planning horizon, in which each decision changes the inventory available to later orders.

\clearpage
\appendix
\makeatletter
\def\@seclabel@section{Appendix~}
\def\@seccntformat#1{\@ifundefined{@seclabel@#1}{}{\@nameuse{@seclabel@#1}}\csname the#1\endcsname\@seccntDot\hskip 0.5em}
\makeatother
\counterwithin{equation}{section}
\counterwithin{table}{section}
\counterwithin{figure}{section}
\section{Forecasting-Method Details}
\label{app:forecast-methods}

This appendix gives the randomized CPS and the Venn--Abers constructions summarized in Section \ref{sec:forecast-methodology}. Throughout, $\boldsymbol x$ is the feature vector of a new forecasting instance.

\paragraph{Randomized CPS.}
With $A_j(\boldsymbol x)$ as in Section \ref{sec:cps} and a random weight $\tau$ drawn uniformly on $[0,1]$, the randomized conformal predictive distribution function counts the shifted outcomes below a candidate deviation $y$ and gives weight $\tau$ to the outcomes equal to $y$ and to the candidate itself:
\begin{equation*}
Q_{\mathrm{CPS}}(y\mid \boldsymbol x;\tau)
=\frac{\left|\{j:A_j(\boldsymbol x)<y\}\right|
+\tau\left(\left|\{j:A_j(\boldsymbol x)=y\}\right|+1\right)}{N_c+1}.
\end{equation*}
Under exchangeability of the calibration and test examples, evaluating $Q_{\mathrm{CPS}}$ at the observed test outcome gives a value uniformly distributed on $[0,1]$ \citep{vovk2020computationally}. The deterministic distribution \eqref{eq:cps-cdf} satisfies $\widehat F(y\mid\boldsymbol x)=[(N_c+1)Q_{\mathrm{CPS}}(y\mid \boldsymbol x;1)-1]/N_c$, which removes the candidate that $Q_{\mathrm{CPS}}$ counts at $\tau=1$.

\paragraph{Inductive Venn--Abers predictor.}
For a binary label and raw model score $\sigma(\boldsymbol x)$, IVAP appends $(\sigma(\boldsymbol x),0)$ to the calibration data and fits a nondecreasing map $f_0$ from model scores to observed class frequencies. It then repeats the calculation with $(\sigma(\boldsymbol x),1)$ to obtain $f_1$. The calibrated probability pair is
\[
[\underline p(\boldsymbol x),\overline p(\boldsymbol x)]=[f_0(\sigma(\boldsymbol x)),f_1(\sigma(\boldsymbol x))].
\]

\paragraph{Cross Venn--Abers predictor.}
Let $\underline{\boldsymbol p}^{(1:n_{\mathrm{cv}})}(\boldsymbol x)$ and $\overline{\boldsymbol p}^{(1:n_{\mathrm{cv}})}(\boldsymbol x)$ collect the calibrated probabilities from the $n_{\mathrm{cv}}$ folds. With $\operatorname{GM}$ the geometric mean and $1-\underline{\boldsymbol p}^{(1:n_{\mathrm{cv}})}(\boldsymbol x)$ taken componentwise, the merged binary forecast is
\begin{equation*}
\widehat p_{\mathrm{VA}}(\boldsymbol x)=
\frac{\operatorname{GM}(\overline{\boldsymbol p}^{(1:n_{\mathrm{cv}})}(\boldsymbol x))}
{\operatorname{GM}(1-\underline{\boldsymbol p}^{(1:n_{\mathrm{cv}})}(\boldsymbol x))+
\operatorname{GM}(\overline{\boldsymbol p}^{(1:n_{\mathrm{cv}})}(\boldsymbol x))}.
\end{equation*}
For the multiclass outcome, let $\widetilde p_k(\boldsymbol x)$ be the merged one-versus-rest probability of class $k$. The normalized class probabilities supplied to the predictive CDF are
\begin{equation*}
\widehat p_k(\boldsymbol x)=
\frac{\widetilde p_k(\boldsymbol x)}{\sum_{k'=1}^{K}\widetilde p_{k'}(\boldsymbol x)},
\qquad k=1,\ldots,K.
\end{equation*}

\section{Proprietary-Data Details}
\label{app:proprietary-features}

This appendix records the feature construction, metric definitions, and secondary diagnostics for the proprietary-data study. The models trained on the JD.com dataset follow the separate setup described in Section \ref{sec:jd-simulator}.

\subsection{Feature Construction and Inputs}

The feature set combines order, stock-keeping-unit (SKU), fulfillment-location, and carrier information. For fulfillment locations, congestion is represented by the estimated waiting time when the order is placed and by shipping volume during the preceding two hours, and recurring workload patterns are summarized within segments defined by wave schedules, operating hours, and facility type. All historical processing-time summaries use training observations only.

For each five-digit ZIP-code (ZIP5) origin--destination pair, carrier features summarize shipment volume, intermediate stops, tracking updates, and carrier-, weather-, and other delay exceptions. The most frequently used origin and destination hubs define a representative route. Spatial features use hub coordinates and Haversine distance. Temporal features summarize the intervals between the main tracking stages. Unseen ZIP5 pairs are imputed from three-digit (ZIP3) and then two-digit (ZIP2) aggregations. Table \ref{tab:all_features} lists the complete inputs.

\begin{table}[!ht]
\centering
\scriptsize
\caption{Forecasting Features Used for the Proprietary Dataset.}
\label{tab:all_features}
\begin{adjustbox}{width=\textwidth}
\begin{tabular}{l | l | l}
\toprule
\textbf{Feature Name} & \textbf{Description} & \textbf{Type} \\ 
\midrule
\multicolumn{3}{c}{\textbf{Order-Level Features}} \\
\midrule
Order Create Month/Day/Hour & Date-time components when order was created & Numerical \\
Order Create Season/Holiday & Season and holiday indicator of order date & Categorical/Binary \\
Promised Delivery Month/Day & Date-time components of promised delivery & Categorical \\
Promised Delivery Holiday & Whether promised date is a holiday & Binary \\
Promised Delivery Time & Days from order to promised delivery & Numerical \\
Shipping Offer Method & Fulfillment experience at order time & Categorical \\
Order Total Units, Weight, Dim. Weight & Size and quantity metrics of the order & Numerical \\
SKU Identifier/Units/Height/Length/Width/Weight & SKU attributes and quantity in order & Mixed \\
SKU Dept/Class/Subclass/Perf. Code & SKU category and performance level & Categorical \\
Customer Membership Tier & Customer membership type & Categorical \\
Dest. ZIP/State & Customer location & Categorical \\
\midrule
\multicolumn{3}{c}{\textbf{Fulfillment Location (FL)-Level Features}} \\
\midrule
FL Identifier/Type/Lat/Lon & Fulfillment location identifier and position & Categorical/Numerical \\
FL to Dest. Distance & Distance to customer & Numerical \\
SKU On-Hand Inventory & SKU availability at location & Numerical \\
FL Capacity Status & Whether FL is over capacity & Binary \\
FL Waiting Units / Same SKU & Queue length for all and same SKU & Numerical \\
FL Shipped Units Last 2Hrs / Same SKU & Recent shipping activity & Numerical \\
FL Processing Window & Time window of order creation & Categorical \\
FL Avg/Std Processing Time & Processing time stats & Numerical \\
\midrule
\multicolumn{3}{c}{\textbf{Carrier-Level Features}} \\
\midrule
Transit Time / Zone / Cutoff Time & Carrier time and region details & Numerical \\
\# Non Transit Days / Before Cutoff & Service schedule and timing compliance & Mixed \\
Route Use Freq. / Median \# Stops/Updates & Route-level shipment frequency & Numerical \\
Total Carrier/Weather/Other Exceptions & Exception event counts & Numerical \\
Avg/Std: Label to Hub, Hub Induction, Hub Transit & Timings between key shipping stages & Numerical \\
Avg/Std: Dest. Scan to Delivery, Final Delivery & Delivery segment durations & Numerical \\
Most Freq. Origin/Dest. Hub & Common route endpoints & Categorical \\
Most Freq. Hub Lat/Lon, Distance & Geolocation and distance between hubs & Numerical \\
\bottomrule
\end{tabular}
\end{adjustbox}
\end{table}

\subsection{Forecast Evaluation Metrics}

Let $N_{\mathrm{test}}$ be the number of test instances, and for instance $i$ let $\boldsymbol x_i$ be its feature vector, $y_i$ its realized delivery-time deviation, and $\widehat F_i(y)=\widehat F(y\mid \boldsymbol x_i)$ its predictive CDF. Its $q$-th quantile is $\widehat y_i^q$ for $q\in(0,1)$, and $\widehat y_i$ is the point prediction being evaluated. Table \ref{tab:evaluation_metrics} gives the formulas used in the forecast study.

\begin{table*}[!ht]
\footnotesize
\centering
\caption{Forecast Evaluation Metrics.}
\label{tab:evaluation_metrics}
\begin{adjustbox}{width=\textwidth}
\begin{tabular}{l|l|l}
\toprule
\textbf{Metric} & \textbf{Formula} & \textbf{Description} \\ \hline
CRPS & $ \frac{1}{N_{\mathrm{test}}}\sum_{i=1}^{N_{\mathrm{test}}} \int_{-\infty}^{\infty} (\widehat F_i(y) - \mathbf{1}\{y \ge y_i\})^2\,dy$ & Continuous ranked probability score.\\ \hline
PL & $\displaystyle
\begin{aligned}
\frac{1}{9}\sum_{q \in \{0.1,0.2,\ldots,0.9\}} &\frac{1}{N_{\mathrm{test}}}\sum_{i=1}^{N_{\mathrm{test}}} \Bigl(q\max\{y_i-\widehat y_i^q,0\} \\
&+(1-q)\max\{\widehat y_i^q-y_i,0\}\Bigr)
\end{aligned}$ & Average pinball loss. \\ \hline
MQCE & $\frac{1}{9}\sum_{q \in \{0.1,0.2,\ldots,0.9\}} \left|\frac{1}{N_{\mathrm{test}}}\sum_{i=1}^{N_{\mathrm{test}}}\mathbf{1}\{y_i<\widehat y_i^q\}-q\right|$ & Mean quantile coverage error. \\ \hline
Accuracy & $\frac{1}{N_{\mathrm{test}}}\sum_{i=1}^{N_{\mathrm{test}}}\mathbf{1}\{y_i=\widehat y_i\}$ & Exact-day accuracy. \\ \hline
RMSE & $\sqrt{\frac{1}{N_{\mathrm{test}}}\sum_{i=1}^{N_{\mathrm{test}}}(\widehat y_i-y_i)^2}$ & Root mean squared error in days. \\ \hline
MAE & $\frac{1}{N_{\mathrm{test}}}\sum_{i=1}^{N_{\mathrm{test}}}|\widehat y_i-y_i|$ & Mean absolute error in days. \\ \hline
\end{tabular}
\end{adjustbox}
\end{table*}

\subsection{Supplementary Forecast Diagnostics}
\label{app:supplementary-forecast-results}

\begin{figure}[pos=!ht]
\centering
    \centering
    \includegraphics[width=0.68\textwidth]{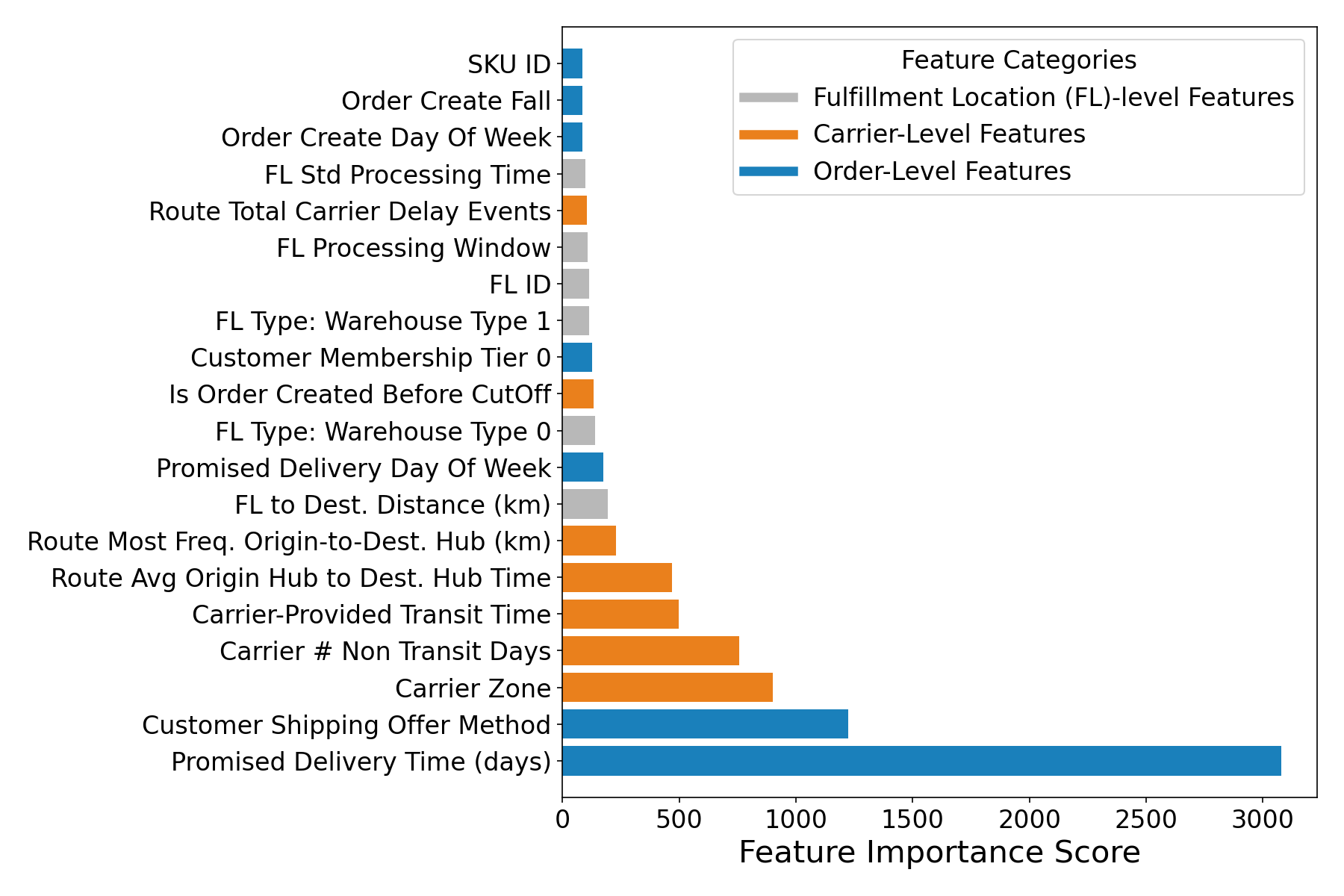}
    \caption{The top 20 features by XGBoost gain importance.}
    \label{fig:feat-imp}
\end{figure}

Figure \ref{fig:feat-imp} ranks the 20 largest XGBoost gain-importance scores. Promised delivery time and the customer shipping-offer method lead the ranking because both describe the delivery promise from which the deviation is measured. Carrier and fulfillment-location variables also appear among the top predictors. These rankings describe how the trained model partitions the data rather than the effect of each feature on delivery time.

\FloatBarrier
\section{Proofs and Supporting Results}
\label{app:theory-proof}

This appendix proves the expected-cost bounds for one decision epoch. Probabilities are taken over calibration, scenario, and candidate-generation randomness.

\subsection{Setup and Method-Specific Bounds}
\label{app:theory-setup}

The analysis compares decisions for the same current order and inventory state. For each decision, the population objective is its expected current-order cost plus expected future-demand recourse cost under the two-stage approximation. The analysis bounds the increase in this objective when C-SAA uses an estimated delivery-time distribution, relative to the best feasible decision under the true delivery-time distribution.

Let $\bar u$ bound the current-order quantity. Assume the feasible decisions and fulfillment options are finite, at least one fulfillment option is available, and $0\leq\sum_{a\in\mathcal A(\boldsymbol s_t)}u_a(\boldsymbol z)\leq\bar u$, with $u_a(\boldsymbol z)\geq0$. Let $\mu_{a,\boldsymbol s_t}$ and $\widehat\mu_{a,\boldsymbol s_t}$ be the true and forecast expected per-unit delivery costs at fulfillment option $a$, each conditional on its context $\boldsymbol x_{a,t}$. All fulfillment-option quantities below are conditional on the same context. The true and forecast population objectives, $J$ and $\widehat J$, share the first-stage order cost and the expected future-demand recourse of \eqref{eq:population-objective}, with future demand following the demand forecast of C-SAA, and differ only in the expected delivery costs:
\begin{align*}
J(\boldsymbol z;\boldsymbol s_t)&=R_t(\boldsymbol z;\boldsymbol s_t)+\E[Q_t(\boldsymbol z,\boldsymbol D_t;\boldsymbol s_t)]+\sum_{a\in\mathcal A(\boldsymbol s_t)}u_a(\boldsymbol z)\mu_{a,\boldsymbol s_t},\\
\widehat J(\boldsymbol z;\boldsymbol s_t)&=R_t(\boldsymbol z;\boldsymbol s_t)+\E[Q_t(\boldsymbol z,\boldsymbol D_t;\boldsymbol s_t)]+\sum_{a\in\mathcal A(\boldsymbol s_t)}u_a(\boldsymbol z)\widehat\mu_{a,\boldsymbol s_t}.
\end{align*}
The comparison assumes that the true delivery distribution at a state is unchanged by additionally conditioning on the offline calibration data and candidate-generation randomness, so $\mu_{a,\boldsymbol s_t}$ does not depend on how the forecasts and candidates were produced.

Let $\mathcal C(\boldsymbol s_t)\subseteq\mathcal Z(\boldsymbol s_t)$ be the nonempty retained candidate set. Let $\boldsymbol z^*$ be a minimizer of $J$ over $\mathcal Z(\boldsymbol s_t)$ and $\boldsymbol z_{\mathcal C}^*$ a minimizer over $\mathcal C(\boldsymbol s_t)$. The candidate-generation gap is
\(
\Gamma_{\mathcal C}^J=J(\boldsymbol z_{\mathcal C}^*;\boldsymbol s_t)-J(\boldsymbol z^*;\boldsymbol s_t).
\)

Both results use the following evaluation conditions. Given the fitted forecasts, calibration data, current state, and candidate-generation and retention decisions, the $N_2$ common evaluation scenarios are independent and identically distributed, with mean objective $\widehat J(\boldsymbol z;\boldsymbol s_t)$ for each candidate. Each scenario objective has conditional range at most $\Delta_{\rm obj}<\infty$. Let $\widehat J_{N_2}$ be the resulting sample average and $\widehat{\boldsymbol z}$ its minimizer over $\mathcal C(\boldsymbol s_t)$. The uniform sampling-error bound $e_{N_2}$ at evaluation failure probability $\alpha_s\in(0,1)$ is
\begin{equation}
e_{N_2}=\Delta_{\rm obj}
\sqrt{\frac{\log(2|\mathcal C(\boldsymbol s_t)|/\alpha_s)}{2N_2}}.
\label{eq:common-sample-error}
\end{equation}

\paragraph{CPS with C-SAA.}

The CPS bound separates uncertainty in the forecast from error in the scenario-based decision procedure. Let $U$ be the scalar response modeled by regression, the deviation in the proprietary study and the duration in the JD study, and let $m$ be the population regression function approximated by the fitted model $h$, so that the population residual is $U-m(\boldsymbol x)$. Let $\mathcal B$ be the family of calibration groups across all fitted forecasting models, with $B=|\mathcal B|$. SCPS has one group per model and MCPS one group per model and prediction range. For group $b$, let $\mathcal J_b$ be the indices of its calibration observations. For fulfillment option $a$, let $b(a)$ be the group used at the current state, $N_c^{b(a)}$ its number of calibration observations, $G_{b(a)}$ the CDF of its population residuals with Lipschitz constant $M_{b(a)}$, and $H_{a,\boldsymbol s_t}$ the CDF of the population residual at the current state. The proper-training information is independent of the calibration and test data and fixes the fitted models and the family of calibration groups. Conditional on this information and on group membership, each group contains $N_c^b\geq1$ calibration observations whose population residuals are iid. The residual empirical-CDF sampling term for calibration failure probability $\alpha_c\in(0,1)$ is $e_{b(a)}=\sqrt{\log(2B/\alpha_c)/(2N_c^{b(a)})}$. The regression errors are $\zeta^{\rm cal}_{b(a)}=\max_{j\in\mathcal J_{b(a)}}|m(\boldsymbol x_j)-h(\boldsymbol x_j)|$ on the calibration observations and $\zeta^{\rm te}_a=|m(\boldsymbol x_{a,t})-h(\boldsymbol x_{a,t})|$ at the current input. The residual mismatch between the calibration group and the current state is $\xi_{a,\boldsymbol s_t}=\sup_{r\in\mathbb R}|G_{b(a)}(r)-H_{a,\boldsymbol s_t}(r)|$. The same map from the response to the delivery-time deviation is applied to true and forecast outcomes, and the resulting cost, as a function of the response, is bounded with finite total variation $V_{a,t}$. The expected-cost error bound for fulfillment option $a$ is
\begin{equation}
\delta^{\rm CPS}_{a,t}=V_{a,t}\bigl[e_{b(a)}+M_{b(a)}\bigl(\zeta^{\rm cal}_{b(a)}+\zeta^{\rm te}_a\bigr)+\xi_{a,\boldsymbol s_t}\bigr].
\label{eq:cps-option-cost}
\end{equation}
The bound applies to the empirical residual sampler defined in Section \ref{sec:cps}.

For SCPS, residual mismatch is measured against the full calibration population, and for MCPS within a fixed group, so grouping can reduce mismatch at the cost of fewer observations per group. The proof controls all fixed groups simultaneously, which allows the group used for an arriving order to depend on earlier decisions.

\paragraph{Venn--Abers with C-SAA.}

For Venn--Abers, the relevant accuracy measure is the error in the probabilities supplied to the optimizer. For this analysis, let $K$ be the number of outcome classes, indexed by $k=1,\ldots,K$. Let $\boldsymbol\theta_a=(\theta_{ak})_{k=1}^K$ be the true class probability vector at the current state and $\boldsymbol p_a=(p_{ak})_{k=1}^K$ the normalized forecast vector, both nonnegative and summing to one. For CVAP, $p_{ak}$ is the normalized probability $\widehat p_k(\boldsymbol x_{a,t})$ of Appendix \ref{app:forecast-methods}. Let $c_{ak}$ be the cost of class $k$ and define its range across classes as $\Delta_a=\max_k c_{ak}-\min_k c_{ak}$.

An excess conditional score measures the expected prediction loss above that attained by the true probabilities at the current state. The excess conditional Brier score is $r^{\rm B}_a=\|\boldsymbol p_a-\boldsymbol\theta_a\|_2^2$, using the sum of squared class errors without division by $K$. The excess conditional log loss is $r^{\log}_a=\sum_k \theta_{ak}\log(\theta_{ak}/p_{ak})$. Terms with $\theta_{ak}=0$ contribute zero. If $p_{ak}=0<\theta_{ak}$, the log loss is infinite. Define the fulfillment-option-level expected-cost error bound $\delta_{a}$ as
\[
\delta_{a}
=\frac{\Delta_a}{2}\min\!\left\{\sqrt{K r^{\rm B}_a},\sqrt{2r^{\log}_a}\right\}.
\]

Either score can be used alone. The bound decreases with the excess Brier score or log loss of the normalized probabilities, and the cost range $\Delta_a$ sets how much a given probability error matters.

\subsection{Auxiliary Stability Lemma}

The following deterministic comparison converts errors in expected delivery costs and candidate evaluation into a bound on the selected decision.

\begin{lemma}[Regret bound with option-specific forecast errors]
\label{lem:forecast-stability}
Let $\delta_a\geq0$ bound the expected-cost error for fulfillment option $a$ and let $\varepsilon\geq0$ bound the candidate-evaluation error. Suppose $|\widehat\mu_{a,\boldsymbol s_t}-\mu_{a,\boldsymbol s_t}|\leq\delta_a$ for every current fulfillment option and
\(
|\widehat J_{N_2}(\boldsymbol z;\boldsymbol s_t)-\widehat J(\boldsymbol z;\boldsymbol s_t)|\leq\varepsilon
\)
for every retained candidate. Then
\[
J(\widehat{\boldsymbol z};\boldsymbol s_t)-J(\boldsymbol z^*;\boldsymbol s_t)
\leq2\bar u\max_{a\in\mathcal A(\boldsymbol s_t)}\delta_a+2\varepsilon+\Gamma_{\mathcal C}^J.
\]
\end{lemma}

\paragraph{Proof.}
For every feasible $\boldsymbol z$, with all sums and maxima over $a\in\mathcal A(\boldsymbol s_t)$,
\[
|\widehat J(\boldsymbol z;\boldsymbol s_t)-J(\boldsymbol z;\boldsymbol s_t)|
\leq\sum_a u_a(\boldsymbol z)|\widehat\mu_{a,\boldsymbol s_t}-\mu_{a,\boldsymbol s_t}|
\leq\bar u\max_a\delta_a.
\]
Because $\widehat{\boldsymbol z}$ minimizes $\widehat J_{N_2}$ on $\mathcal C(\boldsymbol s_t)$, two applications of the common-sample bound give
\[
\widehat J(\widehat{\boldsymbol z};\boldsymbol s_t)
\leq\widehat J(\boldsymbol z_{\mathcal C}^*;\boldsymbol s_t)+2\varepsilon.
\]
Applying the population-objective bound to $\widehat{\boldsymbol z}$ and $\boldsymbol z_{\mathcal C}^*$ and adding
$J(\boldsymbol z_{\mathcal C}^*;\boldsymbol s_t)-J(\boldsymbol z^*;\boldsymbol s_t)=\Gamma_{\mathcal C}^J$
proves the claim. \hfill$\square$

Taking $\delta_a=\delta$ for every fulfillment option in Lemma \ref{lem:forecast-stability} proves Proposition \ref{thm:forecast-decision-loss}.

\subsection{Proof of Proposition \ref{thm:cps-csaa}}

Let $\mathcal B$ be the entire family of calibration groups available to the forecasting models, with $B=|\mathcal B|$. Conditional on independent proper-training information, the fitted models and this family are fixed. Conditional on group-membership indicators, group $b$ contains $N_c^b\geq1$ observations whose population residuals are independent and identically distributed with CDF $G_b$. Conditioning on group membership leaves the calibration covariates random. The bounds below hold simultaneously for all $B$ groups, including groups that the current state does not use.

Let $U$ be the scalar response modeled by regression. For each fitted model, let $m(\boldsymbol x)$ be the population regression function approximated by $h(\boldsymbol x)$. For calibration observation $j\in\mathcal J_b$, let $\boldsymbol X_j$ be its random feature vector and $U_j$ its response. Define the population residual $E_j$, regression error $d_j$, and fitted residual $\widehat E_j$ by
\[
E_j=U_j-m(\boldsymbol X_j),\qquad
d_j=m(\boldsymbol X_j)-h(\boldsymbol X_j),\qquad
\widehat E_j=U_j-h(\boldsymbol X_j)=E_j+d_j.
\]
Assume $G_b$ is $M_b$-Lipschitz. Let $\zeta_b^{\rm cal}=\max_{j\in\mathcal J_b}|d_j|$ and $\zeta_a^{\rm te}=|m(\boldsymbol x_{a,t})-h(\boldsymbol x_{a,t})|$. These errors may be random. The bounds are pointwise in their realized values. Let $b(a)$ denote the group used by fulfillment option $a$ and $H_{a,\boldsymbol s_t}$ the true CDF of $U-m(\boldsymbol x_{a,t})$ at its current state. Define the residual mismatch
\(
\xi_{a,\boldsymbol s_t}=\sup_{r\in\mathbb R}|G_{b(a)}(r)-H_{a,\boldsymbol s_t}(r)|.
\)

Let $\widetilde G_b$ and $\widehat G_b$ be the empirical CDFs of $E_j$ and $\widehat E_j$ over $j\in\mathcal J_b$. Since $|d_j|\leq \zeta_b^{\rm cal}$,
\begin{equation*}
\widetilde G_b(r-\zeta_b^{\rm cal})
\leq\widehat G_b(r)
\leq \widetilde G_b(r+\zeta_b^{\rm cal}),\qquad r\in\mathbb R.
\end{equation*}
Set $e_b=\sqrt{\log(2B/\alpha_c)/(2N_c^b)}$. Conditional on the proper-training information and group memberships, the Dvoretzky--Kiefer--Wolfowitz inequality \citep{dvoretzky1956asymptotic,massart1990tight} gives
\[
\mathbb P\!\left\{
\sup_r|\widetilde G_b(r)-G_b(r)|>e_b
\right\}
\leq2e^{-2N_c^b e_b^2}=\frac{\alpha_c}{B}.
\]
A union bound and then averaging over the conditioning information yield simultaneous coverage for all $B$ groups with probability at least $1-\alpha_c$. On this event, the residual sandwich and Lipschitz property imply
\begin{equation*}
\sup_r|\widehat G_b(r)-G_b(r)|
\leq e_b+M_b \zeta_b^{\rm cal}.
\end{equation*}

Let $\widehat\Phi_{a,t}$ and $\Phi_{a,t}$ denote the forecast and true CDFs of the untransformed response $U$. For a real-valued response threshold $v$, the empirical sampler gives
\(
\widehat\Phi_{a,t}(v)=\widehat G_{b(a)}(v-h(\boldsymbol x_{a,t})),
\)
whereas
\(
\Phi_{a,t}(v)=H_{a,\boldsymbol s_t}(v-m(\boldsymbol x_{a,t})).
\)
Adding and subtracting the shifted versions of $G_{b(a)}$ gives
\begin{equation}
\sup_v|\widehat\Phi_{a,t}(v)-\Phi_{a,t}(v)|
\leq e_{b(a)}+M_{b(a)}
\bigl(\zeta_{b(a)}^{\rm cal}+\zeta_a^{\rm te}\bigr)+\xi_{a,\boldsymbol s_t}.
\label{eq:cps-cdf-bound}
\end{equation}
Because the concentration event covers every group, this inequality also holds when $\boldsymbol s_t$ and its group choices depend on the calibration sample.

To translate response error into cost error, let $\psi_{a,t}$ be a fixed measurable map from the response scale to the delivery-time-deviation scale, applied identically to true and forecast outcomes, and let $g_{a,t}=c_{a,t}\circ\psi_{a,t}$ be the resulting cost as a function of the response. Assume $g_{a,t}$ is bounded with finite total variation $V_{a,t}$. The bounded-variation expectation inequality gives, for any two response distributions $P$ and $P'$ with CDFs $F_P$ and $F_{P'}$,
\[
\left|\int g_{a,t}\,d(P-P')\right|
\leq V_{a,t}\sup_v|F_P(v)-F_{P'}(v)|.
\]
Bounded variation suffices in place of Lipschitz continuity, so $\psi_{a,t}$ may be discontinuous. Applying this inequality to Equation \eqref{eq:cps-cdf-bound} yields
\[
|\widehat\mu_{a,\boldsymbol s_t}-\mu_{a,\boldsymbol s_t}|
\leq V_{a,t}\left[
e_{b(a)}+M_{b(a)}(\zeta_{b(a)}^{\rm cal}+\zeta_a^{\rm te})
+\xi_{a,\boldsymbol s_t}\right]
=\delta^{\rm CPS}_{a,t},
\]
which is the bound in Equation \eqref{eq:cps-option-cost}.

For candidate evaluation, condition on the fitted forecasts, calibration data, the state $\boldsymbol s_t$, and all candidate-generation and retention decisions. Under the stated evaluation assumptions, each candidate's scenario objectives are iid, have mean $\widehat J(\boldsymbol z;\boldsymbol s_t)$, and lie in an interval of width at most $\Delta_{\rm obj}$. With $e_{N_2}$ as in Equation \eqref{eq:common-sample-error} and $\Delta_{\rm obj}>0$, Hoeffding's inequality \citep{hoeffding1963probability} gives, conditionally for each candidate,
\[
\mathbb P\{|\widehat J_{N_2}(\boldsymbol z;\boldsymbol s_t)-\widehat J(\boldsymbol z;\boldsymbol s_t)|>e_{N_2}\}
\leq2\exp\!\left(-\frac{2N_2e_{N_2}^2}{\Delta_{\rm obj}^2}\right).
\]
If $\Delta_{\rm obj}=0$, the evaluation error is zero. A union bound over the retained set, followed by averaging over the conditioning information, bounds the probability of any evaluation error exceeding $e_{N_2}$ by $\alpha_s$. The union bound allows dependent estimates across candidates. Combining this event with the calibration event and applying Lemma \ref{lem:forecast-stability} proves Equation \eqref{eq:cps-costgap}. \hfill$\square$

\paragraph{Scope and asymptotic implication.}
The proposition concerns fixed calibration groups. The empirical MCPS procedure reduces the number of bins when they would be sparse and maps unseen bins to occupied ones. These data-dependent choices lie outside the fixed-group concentration argument and are assessed empirically.

For a sequence of problems indexed by $n$, use a subscript $n$ to denote each problem's quantities and let $\Lambda_n$ be the entire right-hand side of Equation \eqref{eq:cps-costgap}, including all cost-range, regression, residual-mismatch, sampling, and candidate-generation terms. If $\alpha_{c,n}+\alpha_{s,n}\to0$ and $\Lambda_n\to0$ in probability under the joint randomness of calibration and scenario sampling, then the regret converges to zero in probability. Indeed, for every $\epsilon>0$, its exceedance probability is at most $\alpha_{c,n}+\alpha_{s,n}+\mathbb P(\Lambda_n>\epsilon)$. The condition on $\Lambda_n$ includes $\Gamma_{\mathcal C,n}^J\to0$, so candidate generation must also approach the population optimum.

\subsection{Proofs for Venn--Abers with C-SAA}

\paragraph{Proof of Proposition \ref{thm:va-csaa}.}
For a fixed fulfillment-option context, let $\boldsymbol p=(p_k)_{k=1}^K$ and $\boldsymbol\theta=(\theta_k)_{k=1}^K$ be the forecast and true class probability vectors. Let $C$ be the outcome-class index, with $\mathbb P(C=k)=\theta_k$, and write $\mathbb E_{\boldsymbol\theta}$ for expectation under $\boldsymbol\theta$. For any probability vector $\boldsymbol v=(v_k)_{k=1}^K$, define the Brier loss as
\(
L_B(\boldsymbol v,C)=\sum_k(v_k-\mathbf 1\{C=k\})^2.
\)
Expanding the squares gives
\[
\mathbb E_{\boldsymbol\theta}[L_B(\boldsymbol p,C)]-\mathbb E_{\boldsymbol\theta}[L_B(\boldsymbol\theta,C)]
=\sum_k(p_k-\theta_k)^2=\|\boldsymbol p-\boldsymbol\theta\|_2^2.
\]
For log loss, the corresponding excess conditional risk is the Kullback--Leibler (KL) divergence
\(
\sum_k\theta_k\log(\theta_k/p_k)=\operatorname{KL}(\boldsymbol\theta\|\boldsymbol p),
\)
using the zero-probability conventions stated above.

For fulfillment option $a$, use the class costs $c_{ak}$ defined above. Centering these class costs at their midrange and using
$\sum_k(p_{ak}-\theta_{ak})=0$ gives
\begin{equation}
\left|\sum_k(p_{ak}-\theta_{ak})c_{ak}\right|
\leq\frac{\Delta_a}{2}\|\boldsymbol p_a-\boldsymbol\theta_a\|_1.
\label{eq:cost-tv}
\end{equation}
Cauchy--Schwarz and Pinsker's inequality \citep{tsybakov2009introduction} give, respectively,
\[
\|\boldsymbol p_a-\boldsymbol\theta_a\|_1\leq\sqrt{K r_a^{\rm B}},
\qquad
\|\boldsymbol p_a-\boldsymbol\theta_a\|_1\leq\sqrt{2r_a^{\log}}.
\]
The common evaluation event is the event on which every retained candidate's evaluation error is at most $e_{N_2}$. As proved above, it has probability at least $1-\alpha_s$ under the stated evaluation conditions. On this event, substituting either score bound into \eqref{eq:cost-tv} and then into Lemma \ref{lem:forecast-stability} proves Proposition \ref{thm:va-csaa}. \hfill$\square$

The following width bound requires a retained classwise interval before the one-versus-rest probabilities are normalized. For fulfillment option $a$ and class $k$, let $[\underline p_{ak},\overline p_{ak}]$ be that interval and $\widetilde p_{ak}$ the selected probability before normalization. If CVAP uses several folds, define the interval from the retained fold-specific pairs, for example by taking their componentwise envelope. Let $\widetilde{\boldsymbol p}_a=(\widetilde p_{ak})_{k=1}^K$. The lemma treats simultaneous containment of the true conditional class probabilities as an assumption and gives a conditional sensitivity result.

\begin{lemma}[Venn--Abers width sensitivity]
\label{lem:va-width}
Suppose $0\leq\underline p_{ak}\leq\widetilde p_{ak}\leq\overline p_{ak}\leq1$, $\sum_k\widetilde p_{ak}>0$, and, for every class of every current fulfillment option, the conditional probability satisfies $\theta_{ak}\in[\underline p_{ak},\overline p_{ak}]$. Define $w_{ak}=\overline p_{ak}-\underline p_{ak}$ and $W_a^{\rm VA}=\sum_kw_{ak}$. After normalizing $p_{ak}=\widetilde p_{ak}/\widetilde P_a$, where $\widetilde P_a=\sum_{k'}\widetilde p_{ak'}$,
\[
\|\boldsymbol p_a-\boldsymbol\theta_a\|_1\leq2W_a^{\rm VA}.
\]
\end{lemma}

\paragraph{Proof of Lemma \ref{lem:va-width}.}
Simultaneous containment and
$\widetilde p_{ak}\in[\underline p_{ak},\overline p_{ak}]$ imply
\(
|\widetilde p_{ak}-\theta_{ak}|\leq w_{ak}.
\)
Since $\widetilde p_{ak}\geq0$ and $\widetilde P_a>0$,
\begin{align*}
\|\boldsymbol p_a-\boldsymbol\theta_a\|_1
&\leq\|\widetilde{\boldsymbol p}_a/\widetilde P_a-\widetilde{\boldsymbol p}_a\|_1
     +\|\widetilde{\boldsymbol p}_a-\boldsymbol\theta_a\|_1\\
&=|1-\widetilde P_a|+\sum_k|\widetilde p_{ak}-\theta_{ak}|\\
&\leq2\sum_kw_{ak}=2W_a^{\rm VA}.
\end{align*}
\hfill$\square$

By Equation \eqref{eq:cost-tv}, the expected-cost error at each fulfillment option is then at most $\Delta_aW_a^{\rm VA}$, so on the common evaluation event Lemma \ref{lem:forecast-stability} gives the bound \eqref{eq:va-costgap} with $\delta_a$ replaced by $\Delta_aW_a^{\rm VA}$.

\subsection{Marginal Calibration and Conditional Decision Cost}
\label{app:marginal-calibration}

The following example shows that exact marginal calibration does not bound the regret of the resulting decision. Let the context $X\in\{0,1\}$ be equiprobable, let the fulfillment options be $a\in\{0,1\}$, and set $Y_a=\mathbf 1\{a\neq X\}$. Forecast each fixed fulfillment option by a Bernoulli distribution with success probability $1/2$. For each fixed fulfillment option, this forecast matches the true outcome distribution after averaging over $X$, so it is exactly marginally calibrated. Let fulfillment option 0 have zero base cost, fulfillment option 1 have base cost $\eta>0$, and let the late-penalty rate be $\gamma^+>\eta$. The order has one unit, which must be assigned to exactly one fulfillment option, the cost of fulfillment option $a$ is its base cost plus $\gamma^+Y_a$, and decisions are evaluated exactly. The forecast uniquely selects fulfillment option 0. The context-specific optimum selects fulfillment option 0 at $X=0$ and fulfillment option 1 at $X=1$, so the regret, averaged over $X$, is $(\gamma^+-\eta)/2$, which grows without bound with the late-penalty rate while the forecasts remain exactly marginally calibrated.

Because the forecast is constant, it is also calibrated conditional on its own value, the property that Venn--Abers calibration targets, and a constant regression prediction $h\equiv1/2$ places all orders in a single prediction-based Mondrian group. It misstates the expected delivery cost of each fulfillment option by $\gamma^+/2$ in every context, so Proposition \ref{thm:forecast-decision-loss}, with $\bar u=1$ and $\varepsilon=\Gamma_{\mathcal C}^J=0$, bounds the regret at $X=1$ by $\gamma^+$, against an actual regret of $\gamma^+-\eta$. The error that drives the gap depends on the context, which is why Propositions \ref{thm:cps-csaa} and \ref{thm:va-csaa} measure forecast error through the residual mismatch $\xi_{a,\boldsymbol s_t}$ and the conditional scores $r^{\rm B}_a$ and $r^{\log}_a$.

\FloatBarrier
\section{Additional Decision Results}
\label{app:additional-decision-results}

Table \ref{tab:optimization_quantiles} reports the replication quantiles omitted from the headline comparison in Table \ref{tab:optimization_results}. CVAP MC-CLF, MCPS, and SCPS remain below C-SAA with the original nonconformal forecast and both forecast-free benchmarks at each reported percentile.

\begin{table}[!ht]
\centering
\footnotesize
\caption{Replication Quantiles of the Total Six-Day Cost.}
\label{tab:optimization_quantiles}
\begin{tabular}{lrrrr}
\toprule
\textbf{Policy / forecast} & \textbf{Median} $\downarrow$ & \textbf{75th pct.} $\downarrow$ & \textbf{90th pct.} $\downarrow$ & \textbf{95th pct.} $\downarrow$ \\
\midrule
C-SAA + MC-CLF (CVAP) & 800,405 & 808,088 & 816,848 & 819,614 \\
C-SAA + MCPS & 803,530 & 810,090 & 816,314 & 829,527 \\
C-SAA + SCPS & 804,657 & 815,038 & 820,867 & 828,461 \\
C-SAA + MC-CLF (IVAP) & 816,933 & 823,776 & 830,660 & 839,791 \\
C-SAA + original nonconformal forecast & 862,575 & 868,575 & 874,476 & 884,565 \\
Greedy & 925,308 & 936,506 & 951,517 & 959,385 \\
Empirical SAA & 1,006,414 & 1,012,018 & 1,023,105 & 1,029,497 \\
\bottomrule
\end{tabular}
\end{table}

\section{Computational Reproducibility}
\label{app:computational-reproducibility}

This appendix records the remaining forecast and decision implementation details.

\paragraph{Forecast implementation.}
The CPS models use Crepes \citep{bostrom2022crepes}. The CVAP code is adapted from the Venn--Abers Calibration library. The proprietary forecast study uses ten initial MCPS bins. The decision study uses five and reduces the count when a carrier-specific bin would contain fewer than 20 calibration observations. An unseen test bin uses the nearest calibration bin. CVAP calibrates across expanding temporal folds. Calibration labels absent from a fold's proper-training classes are mapped to the nearest observed delivery-day class. Optuna \citep{akiba2019optuna} tunes the base models.

\paragraph{Decision sample.}
The shipping, consolidation, and delivery terms follow \citet{ye2025msom}. The unmet-quantity penalty and sequential inventory dynamics follow \citet{ye2026proxy}. The decision experiment uses March 26--31, 2018. Carrier service 3 has too little data for a carrier-specific model and therefore uses a model covering multiple carrier services, the same fallback for every forecasting method. The other 17 services use carrier-specific forecasts. The delivery penalties are $\gamma^{+}=40.0$ for lateness and $\gamma^{-}=0.2$ for earliness. The lost-sales penalty is $\rho=200$ per unit, and the consolidation discount is 0.5. Shifted CPS duration forecasts are clipped to $[0,5]$ days, and MC-CLF uses the same integer-day duration support. C-SAA rounds each duration draw and subtracts the promised duration before computing delivery penalties. The per-unit shipping cost plus delivery penalty is rounded to three decimal places before multiplication by the assigned quantity.

\paragraph{Candidate generation.}
Following the replicated SAA structure of \citet{kleywegt2002sample} and the implementation of \citet{ye2026proxy}, each candidate solve uses $N_1=50$ paired delivery-cost and future-demand scenarios, and the common evaluation sample has $N_2=500$ scenarios. The implementation targets $L=10$ candidate first-stage decisions from independent $N_1$-scenario samples. It retains candidates for which Gurobi returns an \texttt{OPTIMAL} status within the 30-second time limit, so an order can have fewer than ten retained candidates. Candidate generation and selection therefore use $LN_1+N_2=1{,}000$ sampled scenarios per order. Algorithm \ref{alg:csaa} summarizes the rolling policy, and \citet{ye2026proxy} give the complete first- and second-stage formulation of \eqref{eq:forecast-saa}, including the recourse allocation constraints.

\begin{algorithm}[!ht]
\caption{Rolling C-SAA policy with forecast distribution $\widehat F$.}
\label{alg:csaa}
\begin{algorithmic}[1]
\FOR{each arriving order $t$, with state $\boldsymbol s_t$}
\STATE $\mathcal C(\boldsymbol s_t)\leftarrow\emptyset$
\FOR{$\ell=1,\ldots,L$}
\STATE Draw $N_1$ scenarios of future demand $\boldsymbol D_t^{(\omega)}$ and deviations $Y_{a,t}^{(\omega)}\sim\widehat F(\cdot\mid\boldsymbol x_{a,t})$, $a\in\mathcal A(\boldsymbol s_t)$
\STATE Solve \eqref{eq:forecast-saa} on these scenarios within the time limit
\IF{the solver returns an optimal first-stage decision $\boldsymbol z_\ell$}
\STATE $\mathcal C(\boldsymbol s_t)\leftarrow\mathcal C(\boldsymbol s_t)\cup\{\boldsymbol z_\ell\}$
\ENDIF
\ENDFOR
\STATE Draw a common evaluation sample of $N_2$ scenarios in the same way
\IF{$\mathcal C(\boldsymbol s_t)\neq\emptyset$}
\STATE Implement $\widehat{\boldsymbol z}\in\arg\min_{\boldsymbol z\in\mathcal C(\boldsymbol s_t)}\widehat J_{N_2}(\boldsymbol z;\boldsymbol s_t)$, the objective of \eqref{eq:forecast-saa} averaged over the evaluation sample
\ELSE
\STATE Leave the order unfulfilled
\ENDIF
\STATE Update the state to $\boldsymbol s_{t+1}$
\ENDFOR
\end{algorithmic}
\end{algorithm}

\paragraph{Hardware and solver.}
The decision experiments used a shared compute cluster. The optimization models were solved with Gurobi 13.0.2 through \texttt{gurobipy} 13.0.2, using \texttt{MIPFocus=3}, \texttt{Threads=4}, and a 30-second time limit.

\paragraph{Software environment.}
The experiments used Python 3.12.13, scikit-learn 1.7.2, PyTorch 2.4.1, XGBoost 2.1.4, and CatBoost 1.2.10.
\section*{Acknowledgments}
This research was partly supported by the NSF AI Institute for Advances in Optimization (Award 2112533).

\section*{Declaration of generative AI and AI-assisted technologies in the manuscript preparation process}
During the preparation of this work, the authors used OpenAI Codex and Anthropic Claude Code to assist with manuscript organization, language editing, and LaTeX verification. After using these tools, the authors reviewed and edited the content as needed and take full responsibility for the content of the published article.

\bibliographystyle{cas-model2-names}

\bibliography{refs}

\end{document}